\documentstyle[12pt]{article}

\begin{document}
\title{The non-Archimedean analogs of the 
Bochner-Kolmogorov, Minlos-Sazonov 
and Kakutani theorems.}
\author{Sergey V. Ludkovsky.}
\date{17 October 2000}
\maketitle
\begin{abstract}
Measures on a non-Archimedean Banach space $X$
are considered with values in the real field
$\bf R$ and in the non-Archimedean fields. 
The non-Archimedean
analogs of the Bochner-Kolmogorov and Minlos-Sazonov theorems are given.
Moreover, infinite products of measures are considered and the analog 
of the Kakutani theorem is given.
\end{abstract}
\section{Introduction.}
There are few works about integration in a classical 
Banach space, that is over the field $\bf R$ of real numbers
or the field $\bf C$ of complex numbers
\cite{boui,chri,cons,dal,sko,vah}.
On the other hand, for a non-Archimedean Banach space $X$
(that is over a
non-Archimedean field) this theory is less developed.
An integration in $X$ is a very important part of the 
non-Archimedean analysis.
 The progress of quantum mechanics and different branches of
modern physics related, for example, with theories of elementary
particles lead to the necessity of developing integration theory
in a non-Archimedean Banach space \cite{ish,vla}.
It may also be useful for the development of non-Archimedean analysis.
Non-Archimedean functional analysis develops rapidly
in recent years  and has many principal differences
from the classical functional analysis
\cite{khrum,roo,sch1,sch2,schrn,vla}. 
Topological vector spaces over non-Archimedean fields
are totally disconnected, classes of smoothness
for functions and compact operators are defined
for them quite differently from that of the classical case,
also the notion of the orthogonality of vectors
has obtained quite another meaning.
In the non-Archimedean case analogs of the Radon-Nikodym
theorem and the Lebesgue theorem about convergence
are true under more rigorous and another conditions.
Especially strong differences are for measures 
with values in non-Archimedean fields, because classical
notions of $\sigma $-additivity and quasi-invariance have lost 
their meaning. 
\par On the other hand the development of the 
non-Archimedean functional analysis and its applications in
non-Archimedean quantum 
mechanics \cite{khrum,vla,yojan} leads to the necessity
of solving such problems. For example, problems related with
quantum mechanics on manifolds are related with diffeomorphism groups,
their representations and measures on them \cite{ish,lutmf99}. 
In articles \cite{luumn51,lutmf99,luum983,luseabm} 
quasi-invariant measures on diffeomorphism and loop groups
and also on non-Archimedean manifolds were constructed.
These measures were used for the investigation of irreducible representations 
of topological groups \cite{lutmf99,luseabm,lubp}.
The theorems proved in this work enlarge classes of  
measures on such groups and manifolds,
this also enlarges classes of representations.
For example, theorems of the Minlos-Sazonov type
characterize measures with the help of characteristic functionals
and compact operators. In the non-Archimedean case
compact operators are more useful, than nuclear operators
in the classical case. Theorems of the Bochner-Kolmogorov and
Kakutani type characterize products of measures
and their absolute continuity relative to others measures.
\par In this article measures are considered on Banach spaces,
though the results given below can be developed
for more general topological vector spaces, for example,
it is possible to follow the ideas of works
\cite{macoma,maarma,mamaz}, in which were considered
non-Archimedean analogs of the Minlos-Sazonov theorems for real-valued 
measures on topological vector spaces over non-Archimedean fields
of zero characteristic. But it is impossible to make in one article.
In this article, apart from articles of M\c{a}drecki,
measures are considered also with values in non-Archimedean fields,
for the cases of real-valued measures also Banach spaces
over non-Archimedean fields $\bf K$ of characteristic
$char ({\bf K})>0$ are considered.
It is well-known, that a real-valued measure $m$ on a locally compact 
Hausdorff totally disconnected
Abelian topological group $G$ is called the Haar measure, if 
\par $(H)$ $m(x+A)=m(A)$ for each $x\in G$ and  
each Borel subset $A$ in $G$. \\
For the $s$-free group $G$ a measure $m$ 
with values in a non-Archimedean field $\bf K_s$ satisfy condition $(H)$
only for an algebra of clopen (closed and open) subsets $A$.
Indeed, in the last case if a measure is 
locally finite and $\sigma $-additive on
the Borel algebra of $G$, then it is purely atomic with atoms
being singletons, so it can not be invariant
relative to the entire Borel algebra (see Chapters 7-9 \cite{roo}). 
\par In \S 2 weak distributions, 
characteristic functions of measures and their properties
are defined and investigated.
In \S 3 the non-Archimedean analogs of the Minlos-Sazonov and 
Bochner-Kolmogorov theorems are given. Quasi-measures 
also are considered. Frequently definitions, formulations 
of statements and their proofs differ strongly from that of classical.
In \S 4 products of measures are considered together with their 
density functions. The non-Archimedean analog of the Kakutani theorem is
investigated.
\par  {\bf Notations.} Henceforth, $\bf K$
denotes a locally compact infinite field with a non-trivial norm,
then the Banach space $X$ is over $\bf K$. 
In the present article measures on $X$ have values in $\bf R$
or in the field $\bf K_s$, that is, a finite
algebraic extension of the $s$-adic field $\bf Q_s$ 
with the certain
prime number $s$. We assume that $\bf K$ is $s$-free as the additive group,
for example, either $\bf K$ is a finite algebraic extension of $\bf Q_p$ or 
$char ({\bf K})=p$ and $\bf K$ is isomorphic with a field
$\bf F_p{(\theta )}$ consisting of elements $x=\sum_ja_j\theta ^j$,
where $a_j\in \bf F_p$, $|\theta |=p^{-1}$, 
$\bf F_p$ is a finite field of $p$ elements,
$p$ is a prime number and $p\ne s$. These imply that 
$\bf K$ has the Haar measures with values in $\bf R$ and $\bf K_s$
\cite{roo}. If $X$ is a Hausdorff topological space
with a small inductive dimension $ind (X)=0$, then
\par $E$ denotes an algebra of subsets of $X$, 
as a rule $E\supset Bco(X)$ for $\bf K_s$-valued measures
and $E\supset Bf(X)$ for real-valued measures, where
\par $Bco (X)$ denotes an algebra of clopen (closed and open)
subsets of $X$,
\par $Bf(X)$ is a Borel $\sigma $-field of $X$ in \S 2.1;
\par $Af(X,\mu )$ is the completion of $E$ by a 
measure $\mu $ in \S 2.1;
\par ${\sf M}(X)$ is a space of norm-bounded measures on $X$ in \S 2.1;
\par ${\sf M}_t(X)$ is a space of Radon norm-bounded measures in \S 2.1;
\par $c_0(\alpha ,{\bf K})$ is a Banach space
and $P_L$ is a projector (fixed relative to a chosen basis) in \S 2.2;
\par $\mu _L$ is a projection of a measure $\mu $ in \S 2.2;
\par $\{ \mu _{L(n)}: n \} $ is a sequence of weak distributions in \S 2.2;
\par $B(X,x,r)$ is a ball in \S 2.2;
\par $L(X,\mu ,{\bf K_s})$ in \S 2.4;
\par $\chi _{\xi }$ is a character with values either in $\bf T$ or 
$\bf T_s$ in \S 2.6;
\par $\theta (z)=\hat \mu $ is a characteristic functional in \S 2.6;
\par $\delta _0$ is the Dirac measure in \S 2.8;
\par $\mu _1*\mu _2$ is a convolution of measures in \S 2.11;
\par $\psi _{q,\mu }$ and $\tau _q$ in \S 2.14;
\par $C(X,{\bf K})$ is a space of continuous functions from $X$ 
into $\bf K$ in \S 2.16;
\par $X^*$ is the topological dual space  of $X$ \cite{nari};
\par $\hat {\sf C}(Y,\Gamma )$, $\tau (Y)$ in \S 3.2;
\par $\sf B_{+}$, $\sf C_{+}$ in \S 3.5;
\par $\nu \ll \mu $, $\nu \sim \mu $, $\nu \perp \mu $ in \S 4.1.
\section{Weak distributions and families of measures.}
\par  {\bf 2.1.} For a Hausdorff topological space $X$ with
a small inductive dimension $ind (X)=0$ \cite{eng}
the Borel $\sigma$-field is denoted $Bf(X)$. 
Henceforth, measures $\mu $ are given on a measurable space $(X,E)$.
The completion of $Bf(X)$
relative to $\mu $ is denoted by $Af(X,\mu )$. The total variation of $\mu $
with values in $\bf R$ on a subset $A$ is denoted by 
$\| \mu |_A \| $ or $| \mu |(A)$ for $A \in Af(X,\mu )$.
If $\mu $ is non-negative and $\mu (X)=1$, then it is called a probability
measure.
\par We recall that a mapping $\mu : E \to \bf K_s$ 
for an algebra $E$ of subsets of $X$ is called a measure,
if the following conditions are accomplished:
$$(i)\mbox{ }\mu \mbox{ is additive and } \mu (\emptyset )=0,$$
$$(ii) \mbox{ for each }A \in E\mbox{ there exists the following norm}$$
$\| A\|_{\mu }:=\sup \{|\mu (B)|_{\bf K_s}:
\mbox{  }B \subset A, B\in E\}< \infty $,
$$(iii) \mbox{ if there is a shrinking family } F,
\mbox{ that is, for each }$$
$A, B \in F$ there exist
$F\ni C \subset (A\cap B)$ and $\cap\{A: A \in F\}= \emptyset $, then
$\lim_{A\in F} \mu (A)=0$ (see chapter 7 \cite{roo} and also about the completion
$Af(X, \mu )$ of the algebra $E$ by the measure $\mu $).
A measure with values in $\bf K_s$ is called a probability measure if
$\| X \|_{\mu } =1$ and $\mu (X)=1$. For functions $f: X\to \bf K_s$ and
$\phi : X\to [0, \infty )$
there are used notations $\| f \|_{\phi }:=\sup_{x \in X} (|f(x)|
\phi (x))$,  $N_{\mu }(x):=\inf(\| U\|_{\mu }: \mbox{ }U \in Bco(X),
\mbox{ } x \in X)$, where $Bco(X)$ is an algebra of closed and at
the same time open (clopen)
subsets in $X$. Tight measures (that is, measures defined on
$E\supset Bco(X)$) compose the Banach space ${\sf M}(X)$ with a norm
$\| \mu \|:=\|X\|_{\mu }$.
Everywhere below
there are considered $\sigma $-additive measures with 
$|\mu |(X) < \infty $ and $\|X\|_{\mu }<
\infty $ for $\mu $ with values in $\bf R$ and $\bf K_s$ respectively,
if it is not specified another.
\par A measure $\mu $ on $E$ is called Radon,
if for each $\epsilon >0$ there exists a compact subset $C\subset X$ 
such that $ \| \mu |_{(X\setminus C)} \| < \epsilon $.
Henceforth, ${\sf M}(X)$ denotes a space of norm-bounded
measures, ${\sf M}_t(X)$ is its subspace of Radon
norm-bounded measures.
\par  {\bf 2.2.} Each Banach space $X$ over $\bf K$ 
in view of Theorems 5.13 and 5.16 \cite{roo}
is isomorphic with $c_0(\alpha , {\bf K}):=
\{ x: \mbox{ } x=(x_j:$ $ j \in \alpha ),$
$card(j:$ $|x_j|_{\bf K}>b )< \aleph _0 \mbox{ for each } b>0 \} $, where
$\alpha $ is an ordinal, $card(A)$ denotes the cardinality of $A$,
$\| x \| := \sup (|x_j|: \mbox{ } j \in \alpha )$. A dimension of $X$ over
$\bf K$ is by the definition $dim_{\bf K}X := card (\alpha )$.
For each closed $\bf K$-linear subspace $L$ in $X$ there exists an operator
of a projection $P_L: X\to L$. Moreover, an orthonormal in the
non-Archimedean sense basis in $L$ has a completion to an orthonormal basis
in $X$ such that $P_L$ can be defined in accordance with a chosen basis.
\par If $A \in Bf(X)$, then $P_L^{-1}(A)$ is called a cylinder subset
in $X$ with a base $A$, $B^L:=P_L^{-1}(Bf(L))$,
$B_0:=\cup (B^L:$ $ L\subset X, L  \mbox{ is a Banach subspace ,}
dim_{\bf K}X< \aleph _0 )$.
The minimal $\sigma $-algebra $\sigma B_0$ generated by $B_0$
coincides with $Bf(X)$, if $dim_{\bf K}X\le \aleph _0$. Henceforward, it is
assumed that $\alpha \le \omega _0$, where $\omega _0$ is the initial ordinal
with the cardinality $\aleph _0:=card({\bf N})$. Then there exists
an increasing sequence of Banach subspaces 
$L(n) \subset L(n+1)\subset ...$ such that
$cl(\cup [L(n): n])=X$, $dim_{\bf K}L(n)=\kappa _n $ for each $n$, where
$cl(A)=\bar A$ denotes a closure of $A$ in $X$ for $A\subset X$. 
We fix a family of projections $P^{L(m)}_{L(n)}: L(m)\to L(n)$
such that $P^{L(m)}_{L(n)}P^{L(n)}_{L(k)}=P^{L(m)}_{L(k)}$
for each $m\ge n\ge k$.
A projection
of the measure $\mu $ onto $L$ denoted by
$\mu _L(A):=\mu (P_L^{-1}(A))$ for each $A\in Bf(L)$ compose the consistent
family:
$$(1) \mbox{  } \mu_{L(n)}(A)=\mu _{L(m)}(P_{L(n)}^{-1}(A)\cap L(m))$$
for each 
$m \ge n$, since there are projectors $P_{L(n)}^{L(m)}$,
where $\kappa _n \le \aleph _0$ and there may be chosen
$\kappa _n< \aleph _0$ for each $n$.
\par An arbitrary family of measures $ \{ \mu _{L(n)}: n \in {\bf N} \} $ 
having property $(1)$
is called a sequence of a weak distributions
(see also \cite{dal,sko}).
\par By $B(X,x,r)$ we denote a ball $\{y: \mbox{ }y \in X,
\mbox{ } \|x-y\| \le r \}$, which is clopen (closed and open) in $X$.
\par  {\bf 2.3. Lemma.} {\it A sequence of a weak distributions 
$\{ \mu _{L(n)}: n\}$
is generated by some measure $\mu $ on $Bf(X)$ if and only if
for each $c>0$ there exists $b>0$ such that $| |\mu _{L(n)} |
(B(X,0,r) \cap L(n)) - | \mu _{L(n)} |( L(n))| \le c $ and $\sup_n
|\mu _{L(n)} |(L(n)) < \infty $ for $\mu $ with values in $\bf R$;
\par or $\| L(n) \setminus B(X,0,r) \| _{\mu _{L(n)}} \le c$ and
$\sup_n \|L(n)\|_{\mu _{L(n)}} < \infty $ for $\mu $ with values in
$\bf K_s$, where $r\ge b$.}
\par  {\bf Proof.} In the case of $\mu $ with values in $\bf R$
we can use a Hahn decomposition $\mu =\mu ^+ - \mu ^- $
and substitute everywhere in the proof of Lemma
1 \S 2\cite{sko} a Hilbert space over $\bf R$ onto $X$ over $\bf K$,
since $X$ is a Radon space in view of Theorem 1.2  \S
I.1.3 \cite{dal}, then $|\mu |(A)=\mu ^+(A) +\mu ^- (A)$ for $A \in
Bf(X)$.
\par For $\mu $ with values in $\bf K_s$ the necessity is evident.
To prove the sufficiency it remains only to verify property
(2.1.iii), since then $\|X\|_{\mu }=\sup_n \| L(n)\| _{\mu _{L(n)}}
< \infty $. Let $B(n) \in E(L(n))$, $A(n)=P_{L(n)}^{-1}(B(n))$,
by Theorem 7.6 \cite{roo} for each $c>0$ there is a compact subset
$C(n) \subset B(n)$ such that
$\| B(n)\setminus C(n)\| _{\mu _{L(n)}}<c$, where $\| B(n)\setminus
D(n) \|_{\mu } $ $\le \max(\| B(m)\setminus C(m)\| _{\mu _{L(m)}}: m=1,...,n)
< c$ and $D(n):=\bigcap_{m=1}^n {P_{L(m)}}^{-1}(C(m))\cap L(n))$,
$P_{L(n)}^{-1}(E(L(n))\subset E=E(X)$. If 
$A(n)\supset A(n+1)\supset ...$ and $\bigcap_n A(n)=\emptyset $, then
$A'(n+1)\subset A'(n)$ and $\bigcap_n A'(n)=\emptyset $, where $A'(n):=
P_{L(n)}^{-1}(D(n))$, hence $\| A(n)\| _{\mu } \le \| A'(n)
\| _{\mu }+c$. There
may be taken $B(n)$ as closed subsets in $X$.
In view of the Alaoglu-Bourbaki theorem (see Exer. 9.202(a.3)
\cite{nari}) and the Hahn-Banach theorem (4.8 \cite{roo}) sets
$A(n)$ and $B(X,0,r)$ are weakly compact in $X$, hence, for each $r>0$
there exists $n$ with $B(X,0,r)\cap A(n)=\emptyset $. Therefore,
$\| A(n)\|_{\mu}=\| B(n)\| _{\mu _{L(n)}}$ $\le \| L(n)\setminus
B(X,0,r)\| _{\mu _{L(n)}} \le c$
and there exists $\lim_{n\to \infty } \mu (A(n))=0$,
since $c$ is arbitrary.
\par {\bf 2.4. Definition and notations.} A 
function $\phi :X\to \bf R$ (or $\bf K_s$)
of the form $\phi (x)=\phi _S(P_S x)$ is called a cylinder function if
$\phi _S$ is a $Bf(S)$-measurable (or $E(S)$-measurable respectively)
function on a finite-dimensional over
$\bf K$ space $S$ in $X$. For $\phi _S \in L^1 (S, \mu , {\bf R})$
for $\mu $ with values in $\bf R$ or $\phi _S \in L(S,\mu _S,{\bf K_s}):=
L(\mu _S)$ for $\mu $ with values in $\bf K_s$ we may define an integral 
by a sequence of weak distributions $\{ \mu _{S(n)} \} $: 
$$\int_X \phi (x) \mu _* (dx)
:=\int \phi _{S(n)}(x)\mu _{S(n)} (dx) ,$$ 
where $L(\mu )$ is the Banach space of classes of
$\mu $-integrable functions ($f=g$  $\mu $-almost everywhere, that is,
$\| A\|
_{\mu }=0$, $A:= \{ x:$  $f(x)\ne g(x) \} $ is $\mu $-negligible)
with the following norm $\| f\| :=\| g\| _{N_{\mu }}$ \cite{boui,roo,sko}.
\par {\bf 2.5. Lemma.} {\it A subset $A\subset X=c_0(\omega _0,{\bf K})$
is relatively compact if and only if $A$ is bounded
and for each $c>0$ there exists a finite-dimensional over
$\bf K$ subspace $L\subset X$ such that $\bar A \subset L^c
:= \{ y \in X: \mbox{ } d(y,L):=\inf \{ \| x-y\| : x \in L\} \le c \} $.}
\par {\bf Proof.} If $A$ is bounded and for each $c>0$
there exists $L^c$ with $\bar A \subset L^c$, then there is a
sequence $ \{ k(j): \mbox{ } j \in {\bf N} \} \subset \bf Z$
such that $\lim _{j \to \infty } k(j)= \infty ,$  $\bar A \subset \{
x \in X: \mbox{ } |x_j| \le p^{-k(j)}, \mbox{ } j=1,2,...\}=:S$,
but $X$ is Lindel\"of, $S$ is sequentially compact,
hence $\bar A$ is compact (see \S 3.10.31 \cite{eng}). If
$\bar A$ is compact, then for each $c>0$ there exists a finite number
$m$ such that $\bar A \subset \bigcup_{j=1}^m B(X,x_j,c)$, where
$x_j \in X$. Therefore, $\bar A \subset L^c$ for $L=sp_K(x_j:$
$j=1,...,m):=(x=\sum_{j=1}^m b_jx_j: b_j \in K).$
\par {\bf 2.6. Remarks and definitions.} 
As an additive group $\bf K$ is isomorphic with
$\bf Q_p^n$ with $n \in {\bf N}:= \{ 1,2,... \} $. 
The topologically adjoint space
over $\bf Q_p$
(that is, of continuous linear functionals $f:K\to \bf Q_p$)
is isomorphic with $\bf Q_p^n$ \cite{hew}.
For $x$ and $z \in \bf Q_p^n$ we denote by
$z(x)$ the following sum $\sum_{j=1}^n x_jz_j$, where $x=(x_j:$
$j=1,...,n)$, $x_j \in \bf Q_p$. Each number $y\in \bf Q_p$ has  a
decomposition $y=\sum_l
a_lp^l$, where $\min (l:$ $a_l\ne 0)=:ord_p(y)> - \infty $
($ord(0):=\infty $) \cite{nari}, $a_l \in (0,1,...,p-1)$,
we define a symbol $\{ y \}_p:=\sum_{l<0} a_lp^l$ for
$|y|_p>1$ and $\{ y\} _p=0$ for $|y|_p \le 1$. 
\par For a locally compact field $\bf K$ with
a characteristic $char ({\bf K})=p>0$ let $\pi _j(x):=a_j$
for each $x=\sum_ja_j\theta ^j\in \bf K$ (see Notation).
All continuous characters $\chi : {\bf K}\to \bf C$
(or $\chi : K\to \bf C_s$) have the form 
$\chi =\chi _{\xi } (x)=exp \{ 2\pi i \eta (\xi (x)) \} $,
where $\pi _j: {\bf K}\to \bf R$,
$\eta (x):= \{ x \} _p$ and $\xi \in {\bf Q_p^n}^*=
\bf Q_p^n$ for $char ({\bf K})=0$,
$\eta (x):=\pi _0(x)/p$ and $\xi \in {\bf K}^*=\bf K$
for $char({\bf K})=p>0$,
$x \in \bf K$, $i=(-1)^{(1/2)}$
(see \S 25 \cite{hew}), $exp: {\bf C}
\to \bf C$. 
Each $\chi $ is locally constant, hence $\chi : {\bf K}\to
\bf T$ (or $\chi : {\bf K}\to \bf T_s$) is also continuous, 
where $\bf T$ denotes
the discrete group of all roots of $1$ (by multiplication), $\bf T_s$
denotes its subgroup of elements with orders that are not degrees
$s^m$ of $s$, $m \in \bf N$. 
\par For a measure $\mu $ with values in $\bf R$ or $\bf K_s$
there exists a characteristic functional (that is, called the
Fourier-Stieltjes transformation) $\theta =\theta _{\mu }: C(X,{\bf K}) 
\to \bf C$ or $\bf C_s$:
$$(2) \mbox{  } \theta (f):=\int _X \chi _e(f(x)) \mu (dx),$$
where $e=(1,...,1)$, $x \in X$, $f$ is in the space $C(X,{\bf K})$
of continuous functions from $X$ into $\bf K$, in particular for
$z=f$ in the topologically conjugated space $X^*$
over $\bf K$, $z:X \to \bf K$, $z \in X^*$, $\theta (z)=:
\hat \mu (z)$. 
It has the folowing properties:
$$(3a) \mbox{ }\theta (0)=1 \mbox{ for } \mu (X)=1$$ and $\theta (f)$
is bounded on $C(X,{\bf K})$;
$$(3b)\mbox{ }\sup_f | \theta (f)|=1\mbox{ for probability measures };$$ 
$$(4) \mbox{ } \theta (z) \mbox{ is weakly continuous, that is, } (X^*,
\sigma (X^*,X))\mbox{-continuous},$$
$\sigma (X^*,X)$ denotes a weak topology on $X^*$,
induced by the Banach space $X$ over $\bf K$. 
To each $x \in X$ there corresponds a continuous
linear functional $x^*:X^* \to \bf K$, $x^*(z):=z(x)$, moreover, $\theta
(f) $ is uniformly continuous relative to the norm on 
$$C_b(X,{\bf K}):=\{ f\in C(X,{\bf K}):
\| f\| :=\sup_{x\in X}|f(x)|_{\bf K}<\infty \} ;$$
$$(5) \mbox{ } \theta (z) \mbox{ is positive definite on }X^* \mbox{ and on }
C(X,{\bf K})$$
for $ \mu
\mbox{ with values in } [0, \infty ).$
\par  Property (4) follows from Lemma 2.3, boundedness and continuity of
$\chi _e$ and the fact that due to the Hahn-Banach theorem there is
$x_z \in X$
with $z(x_z)=1$ for $z \ne 0$ such that $z|_{(X \ominus L)}=0$ and
$$\theta (z)=\int_X \chi _e(P_L(x)) \mu(dx)=\int_L \chi _e(y) \mu _L(dy),$$
where $L=Kx_z$, also due to the Lebesgue theorem 2.4.9 \cite{fed}
for real measures (or from Exer. 7.F \cite{roo} for $\mu $ with values in
$\bf K_s$, see also \S 4.2\cite{vah}). Indeed, for each $c>0$ there exists a
compact subset $S \subset X$ such that $|\mu |(X\setminus S)<c$ (or
$\| X\setminus S\| _{\mu }<c$),
each bounded subset $A \subset X^*$ is uniformly equicontinuous
on $S$ (see (9.5.4) and Exer. 9.202 \cite{nari}), that is,
$ \{ \chi _e(z(x)): $ $z \in A \} $ is the uniformly equicontinuous
family (by $x \in S$). On the other hand, $\chi _e(f(x))$
is uniformly equicontinuous on a bounded $A\subset C_b(X,{\bf K})$ by
$x\in S$.
\par Property (5) is accomplished, since 
$$\sum_{l,j=1}^N \theta (f_l
-f_j) \alpha _l \bar \alpha _j =\int_X |\sum_{j=1}^N \alpha _j
\chi _e(f_j(x))|^2 \mu (dx) \ge 0,$$ 
particularly, for $f_j=z_j\in X$, where $\bar \alpha _j$ is a complex 
conjugated number to $\alpha _j.$
\par   We call a functional $\theta $ finite-dimensionally concentrated,
if there exists $L \subset X$, $dim_{\bf K}L< \aleph _0 ,$ such that
$\theta |_{(X\setminus L)}=\mu (X)$. For each $c>0$ and $\delta >0$
in view of Theorem I.1.2 \cite{dal} (or Theorem 7.6\cite{roo}) and Lemma
2.5 there exists a finite-dimensional over $\bf K$ subspace $L$ and compact
$S \subset L^{\delta }$ such that 
$\| X\setminus S\| _{\mu } <c$. Let
$\theta ^L (z):=\theta (P_Lz)$. 
\par This definition is correct, since
$L\subset X$, $X$ has the isometrical embedding into $X^*$ as
the normed space associated with the fixed basis of $X$, 
such that functionals $z \in X$ separate points in $X$.
If $z \in L$, then $|\theta (z)-\theta ^L(z)| \le c\times b \times q$,
where $b=\| X\| _{\mu }$, $q$ is independent of $c$ and $b$.
Each characteristic functional $\theta ^L(z)$ is uniformly continuous by
$z \in L$ relative to the norm $\| *\|$ on $L$, since $|\theta ^L(z)
-\theta ^L(y)| $ $\le |\int_{S' \cap L} [\chi _e(z(x))-\chi _e(y(x))]$
$\mu _L(dx)|$ $+|\int_{L\setminus S'} [\chi _e(z(x))-\chi _e(y(x))]$
$\mu _L(dx)|$, where the second term does not exceed $2C'$ for 
$ \| L\setminus S'\| _{\mu _L}<c'$ for
a suitable compact subset $S' \subset X$ and $\chi _e(z(x))$ is 
an uniformly equicontinuous
by $x\in S'$ family relative to $z\in B(L,0,1)$.
\par Therefore,
$$(6) \mbox{  } \theta (z)=\lim _{n \to \infty } \theta _n(z)$$
for each finite-dimensional over $\bf K$ subspace $L$, where $\theta
_n(z)$ is uniformly equicontinuous and finite-dimensionally concentrated
on $L(n) \subset X$, $z \in X$, $cl(\bigcup_nL(n))=X$, $L(n)
\subset L(n+1)$
for every $n$, for each $c>0$ there are $n$ and $q>0$
such that $|\theta (z)-\theta _j(z)| \le cbq$
for $z \in L(j)$ and $j>n$, $q=const >0$ is independent of
$j$, $c$ and
$b$. Let $ \{ e_j:$ $j \in {\bf N} \} $ be the standard orthonormal 
basis in $X$,
$e_j=(0,...,0,1,0,...)$ with $1$ in $j$-th place. Using
countable additivity of $\mu $, local constantness of $\chi _e$,
considering all $z=be_j$ and $b \in \bf K$, we get that $\theta (z)$ on $X$
is non-trivial, whilst $\mu $ is a non-zero measure, since due to Lemma
2.3 $\mu $ is characterized uniquely by $ \{ \mu _{L(n)} \} $. Indeed,
for $\mu $ with values in $\bf R$ a measure $\mu _V$ on $V$,
$dim_{\bf K}V< \aleph _0 ,$ this follows from
the properties of the Fourier transformation $F$ on spaces of 
generalized functions and also on $L^2(V,\mu _V,{\bf C})$ (see \S 7
\cite{vla}), for $\mu $ with values in
$\bf K_s$ this is also true due to Theorem 9.20\cite{hew}, where 
$$F(g)(z):=
\lim _{r \to \infty } \int_{B(V,0,r)} \chi _e(z(x))g(x)m(dx),$$
$z\in V,$ $g \in L(V,\mu _V,{\bf C_s})$, $m$ is the Haar measure on $V$ 
either with values in $\bf R$ or $\bf K_s$ respectively.
Therefore, the mapping $\mu \mapsto \theta _{\mu }$ is injective.
\par  {\bf 2.7. Proposition.} {\it Let $X=K^j$, $j \in \bf N$,
\par (a) $\mu $ and $\nu $ be real probability measures on $X$,
suppose $\nu $ is symmetric. Then $\int_X \hat \mu (x) \nu (dx)$ $=
\int_X \hat \nu(x) \mu (dx)$ $\in \bf R$ and for each
$0<l<1$ is accomplished the following inequality: \\ 
$\mu ([x \in X:$ 
$\hat \nu (x) \le l]) $ $\le \int_X(1-\hat \mu (x)) \nu (dx)/
(1-l)$. 
\par (b). For each real probability measure $\mu $
on $X$ there exists $r> p^3$ such that for each $R>r$ and $t>0$
the following inequality is accomplished: \\
$\mu ([x \in X:$ $\| x \| \ge tR]) \le $ $
c\int_X[1-\hat \mu (y \xi)] \nu (dy),$ \\
where $\nu (dx)=C\times
exp(-|x|^2)m(dx)$, $m$ is the Haar measure 
on $X$ with values in $[0, \infty )$,
$m(B(X,0,1))=1,$ $\nu (X)=1,$ $2>c=const \ge 1$ is independent on
$t,$ $c=c(r)$ is non-increasing whilst $r$ is increasing, $C>0$.}
\par {\bf Proof.} (a). Recall that $\nu $ is symmetric, if
$\nu (B)=\nu (-B)$ for each $B \in Bf(X)$.
Therefore, $\int_X \chi _e(z(x))\nu (dx)$ $=\int_X \chi _e(-z(x))
\nu (dx)$, that is equivalent to $\int _X sin(2\pi \{ z(x)\} _p)
\nu (dx)=0$
or $\hat \nu (z) \in \bf R$. If $0<l<1$, then $\mu ([x \in X:$ $
\hat \nu (x) \le l])=\mu ([x:$ $1-\hat \nu (x) \ge 1-l])$
$\le \int_X(1-\hat \nu (x)) \mu (dx)/(1-l)=$ $\int_X(1-\hat \mu (x))
\nu (dx)/(1-l)$ due to the Fubini theorem.
\par  (b). Let $\nu (dx)=\gamma (x)m(dx),$ where $\gamma (x)=C\times exp(
-|x|^2),$ $C>0,$ $\nu (X)=1.$
Then $F(\gamma )(z)=:\hat \gamma (z) \ge 0,$ and $\hat \gamma (0)=1$
and $\gamma $ is the continuous positive definite function with
$\gamma
(z)\to 0$ whilst $|z| \to \infty .$ In view of (a):
$\mu ([x: \| x\| \ge tR])$ 
$\le \int_X [1-\hat \mu (y \xi )] \nu (dy)/(1-l),$ 
where $|\xi |=1/t,$ $t>0$, $l=l(R)$.
Estimating integrals, we get (b).
\par {\bf 2.8. Lemma.}{\it Let in the notation of Proposition 2.7
$\nu _{\xi }(dx)=\gamma _{\xi }(x)m(dx)$, $\gamma _{\xi }(x)=
C(\xi )exp(-|x \xi |^2),$ $\nu _{\xi }(X)=1,$ $\xi \ne 0,$ then
a measure $\nu _{\xi }$ is weakly converging to the Dirac measure
$\delta _0$ with the support in $0 \in X$ for $|\xi | \to \infty .$}
\par  {\bf Proof.} We have: $C(\xi )^{-1}=
C_q(\xi )^{-1}=\sum _{l \in \bf Z}[p^{lq}-p^{(l-1)q}]
exp(-p^{2l}|\xi |^2) < \infty ,$ where the sum by $l<0$ does not
exceed $1$, $q=jn$, $j=dim_{\bf K}X,$ $n=dim_{\bf Q_p}{\bf K}.$
Here $\bf K$ is considered as the Banach space 
$\bf Q_p^n$ with the following norm
$|*|_p$ equivalent to $|*|_{\bf K}$,
for $x=(x_1,...,x_j) \in X$ with $x_l \in \bf K$ as usually
$|x|_p=\max_{1\le l \le j}
|x_l|_p$, for $y=(y_1,...,y_n) \in \bf K$  with $y_l \in \bf Q_p$:
$|y|_p:=\max_{1\le l \le n} |y_l|_{\bf Q_p}$.
Further, $p^{l+s}
\sum_{x_l \ne 0} exp(2\pi i\sum_{i=l} ^{-s-1}x_i p^{i+s})$ $=
\int_{p^{l+s}}^1 exp(2 \pi i \phi )d\phi$ $+\beta (s),$ where
$s+l<0$,
$ \lim_{s \to -\infty }(\beta (s)p^{-s-l})=0,$ therefore, $sup[
|\hat \gamma _1(z)|_{\bf R} |z|_X:$ $ z \in X, $ $|z| \ge p^3] \le 2.$
Then taking $0\ne \xi \in \bf K$ and carrying out the substitution of variable
for continuous and bounded functions $f:X \to \bf R$ we get
$\lim _{|\xi | \to \infty }$
$\int_Xf(x) \nu _{\xi }(dx)=f(0).$ This means that $\nu _{\xi }$ is weakly
converging to $ \delta _0$ for $|\xi | \to \infty $.
\par {\bf 2.9. Theorem.}{\it Let $\mu _1$ and $\mu _2$ be 
measures in ${\sf M}(X)$ such that $\hat \mu _1(f)=\hat \mu _2(f)$
for each $f \in \Gamma $. Then $\mu _1 =\mu _2$, where $X=c_0(\alpha ,K),$
$\alpha \le \omega _0$, $\Gamma $ is a vector subspace in a space
of continuous functions $f: X\to \bf K$ separating points in $X$.}
\par {\bf Proof.} Let at first $\alpha < \omega _0$, then due to
continuity of the convolution $\gamma _{\xi }*\mu _j$ by $\xi $,
and Proposition 4.5 \S I.4\cite{vah} and Lemma 2.8 we get
$\mu _1=\mu _2,$ since the family $\Gamma $ generates $Bf(X)$.
Now let $\alpha =\omega _0$, $A=\{x \in
X:\mbox{ } (f_1(x),...,f_n(x)) \in S \},$ $\nu _j$ be an image of a measure
$\mu _j$ for a mapping $x \mapsto (f_1(x),...,f_n(x))$, where either
$S \in Bf({\bf K^n})$ or $S\in E({\bf K^n})$, $ f_j \in
X \hookrightarrow X^*$. Then $\hat \nu _1(y)=\hat \mu _1(y_1f_1+
...+y_nf_n)=\hat \mu _2(y_1f_1+...+y_nf_n)=\hat \nu _2(y)$ for each
$y=(y_1,...,y_n) \in K^n$, consequently, $\nu  _1=\nu _2$ on $E$.
Further we can use the Prohorov theorem 3.4 \S 1.3 \cite{vah},
since compositions of $f\in \Gamma $ with continuous functions
$g: {\bf K}\to \bf R$ or $g: {\bf K}\to \bf K_s$
respectively generate a family of real-valued or $\bf K_s$-valued 
functions correspondingly separating points of $X$.
\par  {\bf 2.10. Proposition.} {\it Let $\mu _l$ and $\mu $
be measures in ${\sf M}(X_l)$ and ${\sf M}(X)$ respectively,
where $X_l=
c_0(\alpha _l,K),$ $\alpha _l \le \omega _0,$ $X=\prod_1^n X_l,$
$n \in \bf N.$ Then the condition $\hat \mu (z_1,...,z_n)=\prod_{l=1}^n
\hat \mu _l(z_l)$ for each $(z_1,...,z_n)\in X\hookrightarrow X^*$
is equivalent to $\mu =\prod_{l=1}^n \mu _l$.}
\par {\bf Proof.} Let $\mu =\prod_{l=1}^n \mu _l$, then
$\hat \mu (z_1,...,z_n)=\int _X \chi _e(\sum z_l(x_l))\prod_{l=1}^n
\mu _l (dx_l)$ $=\prod_{l=1}^n \int_{X_l} \chi _e(z_l(x_l))\mu _l(dx_l).$
The reverse statement follows from Theorem 2.9.
\par {\bf 2.11. Proposition.} {\it Let $X$ be a Banach space over $\bf K$;
suppose $\mu $,
$\mu _1$ and $\mu _2$ are probability measures on $X$. Then
the following conditions are equivalent: $\mu $ is the convolution of two 
measures $\mu _j$, $\mu=\mu _1*\mu _2$, and
$\hat \mu (z)=\hat \mu _1
(z)\hat \mu _2(z)$ for each $z \in X$.}
\par {\bf Proof.} Let $\mu =\mu _1 *\mu _2$. This means by the definition
that $\mu $ is the image of the measure $\mu _1 \otimes \mu _2$ for
the mapping $(x_1,x_2) \to x_1 +x_2,$ $x_j \in X,$ consequently,
$\hat \mu (z)=\int_{X\times X} \chi _e(z(x_1+x_2))$ $(\mu _1 \otimes \mu _2)
(d(x_1,x_2))$ $=\prod_{l=1}^2 
\int_X \chi _e(z(x_l))\mu _l(dx_l)$ $=\hat \mu _1
(z)\hat \mu _2(z).$ On the other hand, if $\hat \mu _1\hat \mu _2=\mu ,$
then $\hat \mu =(\mu _1 *\mu _2)^{\wedge }$ and due to Theorem 2.9 above
for real measures, or Theorem 9.20\cite{roo} for measures with values in
$\bf K_s$, we have $\mu =\mu _1*\mu _2.$
\par {\bf 2.12. Corollary.} {\it Let $\nu $ be a probability measure on
$Bf(X)$ and $\mu *\nu =\mu $ for each $\mu $ with values in the same field,
then $\nu =\delta _0$.}
\par {\bf Proof.} If $z_0 \in X \hookrightarrow X^*$
and $\hat \mu (z_0) \ne 0,$  then from $\hat \mu (z_0)\hat \nu (z_0)=
\hat \mu (z_0)$ it follows that $\hat \nu _0(z_0)=1$. From the property
2.6(6) we get that there exists $m \in \bf N$ with $\hat \mu (z)
\ne 0$ for each $z$ with $\| z \| =p^{-m}$, since $\hat \mu (0)=1$.
Then $\hat \nu (z + z_0)=1$, that is, $\hat \nu |_{(B(X,
z_0,p^{-m}))}=1.$ Since $\mu $ are arbitrary we get
$\hat \nu |_X=1$, that is, $\nu =\delta _0$ due to \S 
2.6 and \S 2.9 for $\bf K_s$-valued
measures and real-valued measures.
\par  {\bf 2.13. Corollary.} {\it Let $X$ and $Y$ be a Banach space 
over $\bf K$,
(a) $\mu $ and $\nu $ be probability measures on $X$ and $Y$
respectively, suppose $T: X\to Y$ is a continuous linear operator.
A measure $\nu $ is an image of $\mu $ for $T$ if and only if
$\hat \nu =\hat \mu \circ T^*,$ where $T^*: Y^* \to X^*$ 
is an adjoint operator. (b). A characteristic functional
of a real measure $\mu $ on $Bf(X)$ is real if and only if $\mu $ is
symmetric.}
\par {\bf Proof} follows from \S 2.6 and \S 2.9.
\par {\bf 2.14. Definition.} We say that a real probability measure
$\mu $ on $Bf(X)$ for a Banach space 
$X$ over $\bf K$ and $0<q< \infty $ has a weak
$q$-th order if $\psi _{q,\mu }(z)=\int_X|z(x)|^q\mu (dx) < \infty $
for each $z \in X^*$. The weakest vector topology in $X^*$ relative to which
all $(\psi _{q, \mu }: \mu)$ are continuous is denoted by $\tau _q$.
\par {\bf 2.15. Theorem.} {\it A characteristic functional $\hat \mu $
of a real probability Radon measure $\mu $ on $Bf(X)$ is continuous in
the topology $\tau _q$ for each $q>0.$}
\par  {\bf Proof.} For each $c>0$ there exists a compact $S
\subset X$ such that $\mu (S)>1-c/4$ and 
$$|1-\hat \mu (z)| \le |\int_S
(1-\chi _e(z(x)))\mu (dx)| + |\int_{X\setminus S} (1-\chi _e(z(x)))
\mu (dx)| \le |1-\hat \mu _c(z)|+c/2,$$ 
where $\mu _c(A)= ( \mu (A \cap S)/
\mu (S)$ and $A \in Bf(X)$; further analogously to the proof of
IV.2.3\cite{vah}.
\par  {\bf 2.16. Proposition.} {\it For a completely regular space
$X$ with $ind(X)=0$ the following statements are accomplished:
\par (a) if $(\mu _{\beta })$ is a bounded net of measures 
in ${\sf M}(X)$ that weakly converges to a measure
$\mu $ in ${\sf M}(X)$, then $(
\hat \mu _{\beta } (f))$ converges to $\hat \mu (f)$ for each continuous
$f: X \to \bf K$; if $X$ is separable and metrizable then $(\hat \mu _{\beta })$
converges to $\hat \mu $ uniformly on subsets that are
uniformly equicontinuous in $C(X,{\bf K})$; 
\par (b) if $M$ is a bounded dense family in a ball of the 
space ${\sf M}(X)$ for 
measures in ${\sf M}(X)$, then a family $(\hat \mu: $ $\mu \in M)$
is equicontinuous on a locally $\bf K$-convex space $C(X,{\bf K})$ in a topology
of uniform convergence on compact subsets $S \subset X$.}
\par {\bf Proof.} (a). Functions $exp(2 \pi i \eta ( \{ f(x) \} ))$ 
are continuous and
bounded on $X$, where $\hat \mu (f)=\int_X \chi _e(f(x))\mu (dx)$. Then
(a) follows from the definition of the weak convergence and Proposition
1.3.9\cite{vah}, since $sp_{\bf C} \{ exp (2\pi i \{ f(x) \} _p):$
$f\in C(X,{\bf K}) \} $ is dense in $C(X,{\bf C})$ and $sp_{\bf C_s} \{ exp
(2\pi i \eta (f(x)):$ $f\in C(X,{\bf K} \} $
is dense in $C(X,{\bf C_s})$.
\par (b). For each $c>0$ there exists a compact subset 
$S \subset X$ such that
$| \mu |(S)>|\mu (X)|-c/4$ for real-valued measures or $\| \mu 
|_{(X\setminus S)} \| <c/4$ for $\bf K_s$-valued measures. 
Therefore, for $\mu \in M$ and $f \in C(X,K)$ with
$|f(x)|_K<c<1$ for $x \in S$ we get $| \mu (X)-Re(\hat \mu (f)|=$ $2
|\int_X sin^2(\pi \eta (f(x)))\mu (dx)|<c/2$
for real-valued $\mu $ and $|\mu (X)-{\hat \mu }(f)|=
|\int_X(1-\chi _e(f(x))\mu (dx)|<c/2$
for $\bf K_s$-valued $\mu $, since for $c<1$ and
$x \in S$ we have $sin (\pi \eta (f(x)))=0$. Further analogously to
the proof of Proposition IV.3.1\cite{vah}, since $X$  is the
$T_1$-space and for each point $x$ and each closed subset $S$ in $X$
with $x \notin S$ there is a continuous function $h: X \to B(K,0,1)$
such that $h(x)=0$ and $h(S)=\{ 1\} $.
\par {\bf 2.17. Theorem.} {\it Let $X$ be a Banach space over $\bf K$,
$\eta :\Gamma \to \bf C$ be a continuous positive definite function,
$(\mu
_{\beta })$ be a bounded weakly relatively compact net in the
space ${\sf M}_t(X)$ of Radon norm-bounded measures 
and there exists $\lim _{\beta } \hat \mu _{\beta } (f)=
\gamma (f)$ for each $f \in \Gamma $ 
and uniformly on compact subsets of the completion $\tilde \Gamma $,
where $\Gamma \subset C(X,K)$
is a vector subspace separating points in $X$. Then
$(\mu _{\beta })$ weakly converges to $\mu \in {\sf M}_t(X)$ with
$\hat \mu |_{\Gamma }=\gamma .$}
\par  {\bf Proof} is analogous to the proof of Theorem IV.3.1\cite{vah}
and follows from Theorem 2.9 above and for $\bf K_s$-valued 
measures using the non-Archimedean
Lebesgue convergence theorem (see Ch. 7 \cite{roo}).
\par  {\bf 2.18. Theorem.} {\it (a). A bounded
family of measures in ${\sf M}({\bf K^n})$ 
is weakly relatively compact if and only if a family
$(\hat \mu: $ $\mu \in M)$ is equicontinuous on $\bf K^n$. 
\par (b). If $(\mu _j:$ $j \in {\bf N})$ 
is a bounded sequence of measures
in ${\sf M}_t({\bf K^n})$, $\gamma : {\bf K^n} \to \bf C$
is a continuous (and in addition positive definite 
for real-valued $\mu _j$) function, 
$\hat \mu _j(y) \to \gamma (y)$ for each $y \in \bf K^n$
(and uniformly on compact subsets in $\bf K^n$
for $\bf K_s$-valued measures),
then $(\mu _j)$ weakly converges to a measure
$\mu $ with $\hat \mu =\gamma .$ 
\par (c). A bounded sequence of measures
$(\mu _j)$ in ${\sf M}_t({\bf K^n})$ 
weakly convereges to a measure
$\mu $ in ${\sf M}_t({\bf K^n})$ if and only if for each $y \in K^n$
there exists $\lim _{j \to \infty }\hat \mu _j(y)=\hat \mu (y).$
\par (d). If a bounded net $(\mu _{\beta })$ in ${\sf M}_t({\bf K^n})$
converges uniformly on each bounded subset
in $\bf K^n$, then $(\mu _{\beta })$ converges weakly to a
measure $\mu $ in ${\sf M}_t({\bf K^n})$, where $n \in \bf N$.}
\par {\bf Proof.} (a). This follows from the Prohorov theorem 1.3.6\cite{vah}
and Propositions 2.7, 2.16. 
\par (b). We have the following inequality:
$\lim_m \sup_{j>m}
\mu _j([x \in K^n:$ $|x| \ge
tR]) $ $\le 2 \int_{K^n} (1-Re(\eta (\xi y)))\nu (dy)$ 
with $|\xi |=1/t$ due to
\S 2.7 and \S 2.8 for real-valued measures. 
Due to the non-Archimedean Fourier transform and the Lebesgue 
convergence theorem \cite{roo} for $\bf K_s$-valued measures and
from the condition $\lim_{R\to \infty } \sup_{|y|>R} | \gamma (y)|R^n=0$
it follows, that for each $\epsilon >0$ there exists $R_0>0$ such that 
$\lim_m \sup_{j>m} \| \mu _j |_{ \{ x\in {\bf K^n}: |x|>R \} } \| 
\le 2\sup_{|y|>R}|\gamma (y)|R<\epsilon $ for each $R>R_0$.
In view of Theorem 2.17 $(\mu _j)$ converges weakly to
$\mu $ with
$\hat \mu =\gamma $. (c,d). These may be proved analogously to IV.3.2\cite{vah}.
\par  {\bf 2.19. Corollary.} {\it If $({\hat \mu }_{\beta })\to 1$ 
uniformly on some
neighbourhood of $0$ in $\bf K^n$ for a bounded net
of measures $\mu _{\beta } $ in ${\sf M}_t({\bf K^n})$, then
$(\mu _{\beta })$ converges weakly to $\delta _0$.}
\par  {\bf 2.20. Definition.}  A family of probability measures
$M \subset
{\sf M}_t(X)$ for a Banach space 
$X$ over $\bf K$ is called planely concentrated if
 for each $c>0$ there exists a $\bf K$-linear subspace $S \subset X$ with
$dim_{\bf K}S=n< \aleph _0 $ such that $\inf(\mu
(S^c))|: $ $\mu \in M)>1-c$. The Banach space ${\sf M}_t(X)$ is 
supplied with the following norm
$\| \mu \|:=|\mu |(X).$
\par {\bf 2.21. Lemma.} {\it Let $S$ and $X$ be the same as in \S 2.20;
$z_1,...,z_m \in X^*$ be a separating family of points in $S$. Then a set
$E:=S^c \cap (x \in X:$ $|z_j(x)| \le r_j; $ $j=1,...,m)$ is bounded
for each $c>0$ and $r_1,...,r_m\in (0,\infty ).$}
\par  {\bf Proof.} A space $S$ is isomorphic with $\bf K^n$, consequently,
$p(x)=\max(|z_j|:$ $j=1,...,m)$ is a norm in $S$ equivalent to
the initial norm.
\par {\bf 2.22. Theorem.} {\it Let $X$ be a Banach space 
over $\bf K$ with a family
$\Gamma \subset X$ separating points in $M\subset {\sf M}_t(X)$. Then $M$
is weakly relatively compact if and only if a family $ \{ \mu _z:$
$\mu \in M \} $ is weakly relatively compact for each $z \in \Gamma $
and $M$ is planely concentrated, where $\mu _z$ is an image measure 
on $\bf K$ of a measure $\mu $ induced by $z$.}
\par  {\bf Proof} follows from Lemmas 2.5, 2.21 and the Prohorov theorem
(see also Theorem 1.3.7\cite{vah} with a substitution $[-r_j,r_j]$
onto $B({\bf K},0,r_j)$).
\par {\bf 2.23. Theorem.} {\it For $X$ and $\Gamma $ the same as
in Theorem 2.22 a sequence $ \{ \mu _j:$ $j \in {\bf N} \} \subset
{\sf M}_t(X)$  is weakly convergent to $\mu \in {\sf M}_t(X)$
if and only if for each $z \in \Gamma $ there exists
$\lim _{j\to \infty } {\hat \mu }_j(z)={\hat \mu }(z)$ and a family
$\{ \mu _j \} $ is planely concentrated.}
\par  {\bf Proof}  follows from Theorems 2.17,18,22 (see also
Theorem IV.3.3\cite{vah}).
\par  {\bf 2.24. Proposition.} {\it Let $X$ be a weakly regular
space with $ind(X)=0$, $\Gamma \subset C(X,K)$ be a vector subspace
separating points in $X$, $(\mu _n:$ $ n \in {\bf N})$ $\subset
{\sf M}_t(X),$ $\mu \in {\sf M}_t(X)$, $\lim_{n \to \infty }
\hat \mu _n(f)=\hat \mu (f)$ for each $f \in \Gamma $. Then
$(\mu _n)$ is weakly convergent to $\mu $ relative to the weakest
topology $\sigma (X,\Gamma )$ in $X$ relative to which all
$f \in \Gamma $ are continuous.}
\par  {\bf Proof} follows from Theorem 2.18 and is analogous to
the proof of Proposition IV.3.3\cite{vah}.
\section{The non-Archimedean analogs of the Minlos-Sazonov 
and Bochner-Kolmogorov theorems.}
\par  {\bf 3.1.} Let $(X,{\sf U})$ $=\prod_{\lambda }
(X_{\lambda },{\sf U}_{\lambda })$ be a product of measurable
completely regular Radon spaces $(X_{\lambda },{\sf U}_{\lambda })$
$=(X_{\lambda },{\sf U}_{\lambda },{\sf K}_{\lambda })$, where
${\sf K}_{\lambda }$ are compact classes approximating from below
each measure $\mu _{\lambda }$ on $(X_{\lambda },{\sf U}_{\lambda })$,
that is, for each
$c>0$  and elements $A$ of an algebra ${\sf U}_{\lambda }$
there is $S \in {\sf K}_{\lambda }$, $S \subset A$
with $\| A\setminus S\| _{\mu _{\lambda }} <c$. 
\par {\bf Theorem.} {\it Each
bounded quasi-measure $\mu $ with values in $\bf K_s$ on $(X,{\sf U})$
(that is, $\mu |_{\sf U_{\lambda }}$ is a bounded measure 
for each $\sf \lambda $)
is extendible to a measure on an algebra $Af(X,\mu ) \supset \sf U$,
where an algebra $\sf U$ is generated by a family $({\sf U}_{\lambda }:$
$\lambda \in\Lambda )$.}
\par {\bf Proof.} We have 2.1(i) by the condition and $\| X\| _{\mu }
< \infty $, if 2.1(iii) is satisfied. It remains to prove 2.1(iii).
For each sequence $(A_n) \subset \sf U$ with $\bigcap _n
A_n=\emptyset $ and each $c>0$ for each $j \in \bf N$ we choose
$K_j \in \sf K$, where the compact class $\sf K$ is generated by
$({\sf K}_{\lambda })$ (see Proposition 1.1.8\cite{dal}), such that
$K_j\subset A_j$ and $\| A_j\setminus K_j\| _{\mu } <c$. Since
$\bigcap _{n=1}^{\infty } K_n $ $\subset \bigcap _n A_n=\emptyset $,
then there exists $l \in \bf N$ with $\bigcap _{n=1}^l K_n=\emptyset $,
hence $A_l=
A_l\setminus \bigcap _{n=1}^l K_n$ $\subset \bigcup _{n=1}^l
(A_n\setminus K_n)$, consequently, $\| A_l\| _{\mu } \le \max_{n=1,...,l}
(\| A_n\setminus K_n\| _{\mu })<c$. It remains to use Theorem 7.8\cite{roo}
about uniqueness of an extension of a measure.
\par {\bf 3.2. Definition.} Let $X$ be a Banach space 
over $\bf K$, then a mapping
$f:X\to \bf C$ is called pseudocontinuous, if its restriction $f|_L$
is uniformly continuous for each subspace  $L\subset X$ with
$dim_{\bf K}L<\aleph _0$. Let $\Gamma $ be a family of mappings $f:Y\to \bf K$
of a set $Y$ into a field $\bf K$. We denote by $\hat
{\sf C}(Y,\Gamma )$ the minimal $\sigma $-algebra (that is called
cylinder) generated by an algebra ${\sf C}(Y,\Gamma )$
of subsets of the form $C_{f_1,...,f_n;E}:= \{ x\in X:$ $(f_1(x),...,f_n(x))
\in S \} $, where $S \in Bf({\bf K^n})$, $f_j \in \Gamma $. 
We supply $Y$ with
a topology $\tau (Y)$ which is generated by a base
$(C_{f_1,...,f_n;E}:$ $f_j \in \Gamma ,$
$ E $ $\mbox{is open in}$ $K^n)$.
\par  {\bf 3.3. Theorem. Non-Archimedean analog of the Bochner-Kolmogorov
theorem.} {\it Let $X$ be a Banach space 
over $\bf K$, $X^a$ be its algebraically
adjoint $\bf K$-linear space (that is, of all linear mappings
$f:X\to \bf K$ not necessarily continuous). A mapping $\theta : X^a\to \bf C$
is a characteristic functional of a probability measure
$\mu $ with values in $\bf R$ [or $\bf K_s$] and is defined on
$\hat {\sf C}(X^a,X)$ [or $C(X^a,X)$]
if and only if $\theta $ satisfies conditions
2.6(3,5) for $(X^a,\tau (X^a))$ and is pseudocontinuous on $X^a$ [
or $\theta $ satisfies 2.6(3,6) for $(X^a,\tau (X^a)$ and is
pseudocontinuous on $X^a$
respectively].}
\par  {\bf Proof.} (I). For $dim_{\bf K}X=card(\alpha )<\aleph _0$
a space $X^a$ is isomorphic with $K^{\alpha }$, hence the statement of
theorem for a measure $\mu $ with values in $\bf K_s$ follows from Theorem
9.20\cite{roo} and Theorems 2.9 and 2.18 above, since
$\theta (0)=1$ and $|\theta (z)|\le 1$ for each $z$.
\par (II). We consider now the case of $\mu $ with values in $\bf R$
and $\alpha < \omega _0$. In \S 2.6 
(see also \S 2.16-18,24) it was proved that $\theta =\hat \mu $
has the desired properties for real probability measures $\mu $.
On the other hand, there is $\theta $ which satisfies the conditions of
the theorem. Let $\theta _{\xi } (y)=\theta (y)h_{\xi }(y)$, where
$h_{\xi }(y)=F[C(\xi )exp(-\| x \xi \| ^2)](y)$ (that is, the Fourier transform
by $x$), $\nu _{\xi }({\bf K^{\alpha }})=1$, $\nu _{\xi }(dx)=C(\xi )
exp(-\| x \xi \| ^2)m(dx)$ (see Lemma 2.8), $\xi \ne 0$.
Then $\theta _{\xi }(y)$ is positive definite and is uniformly continuous
as a product of two such functions. Moreover, $\theta
_{\xi }(y) \in L^1({\bf K^{\alpha }},m,{\bf C})$. For $\xi \ne 0$
a function $f_{\xi }(x)=$ $\int_{\bf K^{\alpha }} \theta _{\xi }(y)
\chi _e(x(y))m(dy)$ is bounded and continuous, a function
$exp(-\| x \xi \| ^2)=:s(x)$ is positive definite. Since $\nu _{\xi }$
is symmetric and weakly converges to $\delta _0$,
hence there exists $r>0$ such that for each $|\xi |>r$ we have
${\hat \gamma }_{\xi }(y)=$ $\int _{\bf K^{\alpha }} C(\xi )exp(-\| $ $x \xi
\| _p^2)exp(2 \pi i \eta y(x)))m(dx)$ $=\int cos(2\pi \eta (y(x)))
exp(-\| x\eta \| _p^2)C($ $\xi )m(dx)/2$ $>1-1/R$ for $|y| \le R$,
consequently, ${\hat \gamma }_{\xi }(y)={\hat \zeta }_{\xi }^2(y)$
for $|y| \le R$, where ${\hat \zeta }_{\xi }$ is positive definite
uniformly continuous and has a uniformly continuous extension on
$\bf K^{\alpha }$. Therefore, for each $c>0$ there exists $r>0$ such that
$\| \nu _{\xi }-\kappa _{\xi }*\kappa _{\xi } \| <c$ for each
$|\xi |>r$, where $\kappa _{\xi }(dx)=\zeta _{\xi }(x)m(dx)$
is a $\sigma $-additive non-negative measure.
Hence due to corollary from Proposition IV.1.3\cite{vah}
there exists $r>0$
such that $\int_{\bf K^{\alpha }} \theta _{\xi }(y) \chi _e(-x(y))
\nu _j(dy) \ge 0$ for each $|j|>r$, consequently,
$f_{\xi }(x)=\lim _{|j| \to \infty } \int _{\bf K^{\alpha }}
\theta _{\xi }(y)$ $ \chi _e(-x(y))\nu _j(dy) \ge 0$.
From the equality $F[F(\gamma _{\xi })(-y)](x)=\gamma _{\xi }(x)$
and the Fubini theorem it follows that $\int f_{\xi }\chi _e(y(x))h_j
(x)m(dx)$ $=\int
\theta _{\xi }(u+y)\nu _j(du)$. For $y=0$ we get
$\lim _{|\xi | \to \infty } \int f_{\xi }(x)m(dx)=$ $\int f(x)m(dx)$
$=\lim_{|\xi | \to \infty } \lim_{|j| \to \infty } \int f_{\xi }
(x)h_j(x)m(dx)$ and  $\lim_{|\xi |} \lim_{|j|} |\int_{\bf K^{\alpha }}
\theta _{\xi }(u)\nu _j(du)| \le 1$. From Lemma 2.8 it follows that
$\hat f(y)=\theta (y)$, since by Theorem 2.18 $\theta =\lim _{|\xi |
\to \infty } \theta _{\xi }$ is a characteristic function of a probability
measure on $Bf({\bf K^{\alpha }})$, where $f(x)=
\int_{\bf K^{\alpha }} \theta (y)\chi _e(-x(y))m(dy)$.
\par  (III). Now let $\alpha =\omega _0$. It remains to show that
the conditions imposed on $\theta $ are sufficient, because their necessity
follows from the modification of 2.6 (since $X$ has an algebraic embedding
into $X^a$).
The space $X^a$ is isomorphic with ${\bf K}^{\Lambda }$ which is the space
of all $\bf K$-valued functions defined on the Hamel basis $\Lambda $
in $X$. The Hamel basis exists due to the Kuratowski-Zorn lemma
(that is, each finite system of vectors in $\Lambda $ is linearly independent
over $\bf K$, each vector in $X$ is a finite linear combination over
$\bf K$ of elements from $\Lambda $). Let $J$ be a family of all non-void
subsets in $\Lambda $. For each $A \in J$ there exists a functional
$\theta _A: {\bf K}^A \to \bf C$ such that $\theta _A(t)=
\theta (\sum_{y \in A}t(y)y)$ for $t \in {\bf K}^A$. From the conditions
imposed on $\theta $ it follows that $\theta _A(0)=1$, $\theta _A$ is
uniformly continuous and bounded on ${\bf K}^A$, 
moreover, it is positive definite
(or due to 2.6(6) for each $c>0$ there are $n$ and $q>0$ such that
for each $j>n$ and $z \in {\bf K}^A$ the following inequality is satisfied:
$$(7)\mbox{  }|\theta _A(z)-\theta _j(z)| \le cbq,$$
moreover, $L(j) \supset {\bf K}^A$, $q$ is independent on $j$, $c$ and $b$.
From (I,II) it follows that on $Bf({\bf K}^A)$ 
there exists a probability measure
$\mu _A$ such that $\hat \mu _A =\theta _A$. The family of measures
$ \{ \mu _A:$ $A \in J \} $ is consistent and bounded, since
$\mu _A=\mu _E \circ
(P_E^A)^{-1}$, if $A \subset E$, where $P_E^A: {\bf K}^E \to
{\bf K}^A$ are the natural projectors. Indeed, in the case of measures with
values in $\bf R$ each $\mu _A$ is the probability measure. For
measures with values in $\bf K_s$ this is accomplished due to conditions
(7), 2.6(6) for $X^a$ and due to Theorem 9.20\cite{roo}.
\par  In view of Theorem 1.1.4\cite{dal} (or Theorem 3.1 above) on
a cylinder $\sigma $-algebra of the space ${\bf K}^{\Lambda }$
there exists the unique measure $\mu $ such that $\mu _A=\mu
\circ (P^A)^{-1}$ for each $A \in J$, where $P^A: {\bf K}^{\Lambda }\to 
{\bf K}^A$
are the natural projectors. From $X^a={\bf K}^{\Lambda }$ it follows that
$\mu $ is defined on $\hat {\sf C}(X^a,X)$
(or on $C(X^a,X)$ for $\bf K_s$-valued measures).
For $\mu $ on ${\hat C}(X^a,X)$ or $C(X^a,X)$ there exists its extension
on $Af(X,\mu )$ such that $Af(X,\mu )\supset Bco(X)$ (see \S 2.1).
\par  {\bf 3.4. Definition. \cite{sch2} } A continuous linear operator
$T:X\to Y$ for Banach spaces
$X$ and $Y$ over $\bf K$ is called compact, if
$T(B(X,0,1))=:S$ is a compactoid, that is, for each neighbourhood
$U \ni 0$ in $Y$ there exists a finite subset
$A \subset Y$ such that $S \subset U
+co(A)$, where $co(A)$ is the least $\bf K$-absolutely convex subset
in $V$ containing $A$ (that is,
for each $a$ and $b\in \bf K$ with $|a|\le 1$,
$|b| \le 1$ and for each $x,y \in V$ the following inclusion
$ax+by \in V$ is accomplished).
\par  {\bf 3.5.} Let $\sf B_+$ be a subset of non-negative
functions which are $Bf(X)$-measurable and let $\sf C_+$ be its subset
of non-negative cylinder functions. By ${\hat {\sf B}}_+$ we denote
a family of functions $f \in \sf B_+$ such that $f(x)=\lim_{n}g_n(x)$,
$g_n\in \sf C_+$, $g_n \ge f$. For $f \in {\hat {\sf B}}_+$ let
$\int_X f(x)\mu _*(dx)=\inf_{g \ge f, \mbox{  }g \in {\sf C_+}}
\int_Xg(x)\mu _*(dx)$.
\par For $f\in L(X,\mu ,{\bf K_s})$ and $\bf K_s$-valued measure
$\mu $ let $\int_Xf(x)\mu _*(dx)=\lim_{n\to \infty }\int_Xg_n(x)\mu _*(dx)$
for norm-bounded sequence of cylinder functions
$g_n$ from $L(X,\mu ,{\bf K_s})$
converging to $f$ uniformly on compact subsets of $X$.
Due to the Lebesgue converging theorem this limit exists and does not depend
on a choice of $\{ g_n: n \} $.
\par  {\bf 3.6. Lemma.} {\it A sequence of a weak distributions
$(\mu _{L(n)})$
of probability Radon measures is generated by a real probability neasure
$\mu $ on $Bf(X)$ of a Banach space 
$X$ over $\bf K$ if and only if there exists
$$(8)\mbox{  }\lim_{|\xi |\to \infty }\int_X G_{\xi }(x)\mu_*(dx)=1,$$
where $\int_X G_{\xi }(x)\mu _*(dx):=S_{\xi }(\{ \mu _{L(n)}:n \})$  and \\
$S_{\xi }(\{ \mu _{L(n)} \}):=\lim_{n \to \infty }$ $\int_{L(n)}
F_n(\gamma _{\xi ,n})(x)$ $\mu _{L(n)}(dx),$
$\gamma _{\xi ,n}(y):=\prod_{l=1}^{m(n)} \gamma _{\xi }(y_l)$, \\
$F_n$ is a Fourier transformation by $(y_1,...,y_n)$, $y=(y_j:$ $j \in \bf N)$,
$y_j \in \bf K$,
$\gamma _{\xi }(y_l)$ are the same as in Lemma 2.8 for ${\bf K}^1$; here
$m(n)=dim_{\bf K}L(n)< \aleph _0$, $cl(\bigcup_n L(n)=X=c_0(\omega _0,K)$.}
\par  {\bf Proof.} If a sequence of weak ditributions 
is generated by a measure
$\mu $, then in view of 2.6(3-6), Lemmas 2.3, 2.5, 2.8, Propositions
2.10 and 2.16, Corollary 2.13, the Lebesgue convergence theorem
and the Fubini theorem, also from the proof of Theorem 3.3 and the Radon
property of $\mu $
it follows that there exists $r>0$ such that 
$$\int_X G_{\xi }
(x)\mu _*(dx)=\int_X G_{\xi }(x)\mu (dx)=\lim_{n \to \infty }
\int_{L(n)} \gamma _{\xi ,n}(y) \hat \mu _{L(n)}(y) m_{L(n)}(dy),$$
since $\lim _{j \to \infty }
x_j=0$ for each $x=(x_j:$ $ j) \in X$. In addition, $\lim _{|\xi |
\to \infty }
S_{\xi }(\{ \mu _{L(n)} \})$ $=\int_X \mu (dx)=1$. Indeed,
for each $c>0$ and $d>0$ there exists a compact $V_c \subset X$
with $\| \mu |_{(X\setminus
V_c)} \| <c$ and there exists $n_0$ with $V_c \subset L(n)^d$ for each
$n>n_0$.
Therefore, choosing suitable sequences of $c(n)$, $d(n)$, $V_{c(n)}$ and
$L(j_n)$ we get that $[\int_{L(n)} \gamma _{\xi ,n}(y) \hat \mu _{L(n)}
(y) m_{L(n)}(dy):$ $n \in {\bf N}]$ is a Cauchy sequence, where
$m_{L(n)}$ is the real Haar measure on $L(n)$, the latter is considered as
$\bf Q_p^{m(n)b}$, $b=dim_{\bf Q_p}{\bf K}$, $m(B(L(n),0,1)=1$.
Here we use $G_{\xi }(x)$ for a formal expression of the limit
$S_{\xi }$ as the integral. Then $G_{\xi }(x)$
$(mod$ $\mu )$ is defined evidently as a function for $\mu $ or
$ \{ \mu_{L(n)}:$ $n \} $ with a compact support, also for  $\mu $ with
a support in a finite-dimensional subspace $L$ over $\bf K$ in $X$.
By the definition $supp(\mu _{L(n)}:$ $n)$ is compact, if
there is a compact $V
\subset X$ with $supp(\mu _{L(n)}) \subset P_{L(n)}V$ for each
$n$.
That is, condition (8) is necessary.
\par On the other hand, if (8) is satisfied, then for each
$c>0$ there exists $r>0$ such that $|\int_X G_{\xi }(x) \mu _*(dx)-1|<
c/2$ for real-valued measures or $ | \| G_{\xi }(x) \| _{\mu _*}-1 |<c/2$
for $\bf K_s$-valued-measures, when
$|\xi | >r$, consequently, there exists $n_0$ such that for each
$n>n_0$ the following inequality is satisfied:
$$|1-\int_XF_n(\gamma _{\xi ,n})(x)\mu _*(dx)|
\le | \| \mu |_{(L(n)\cap B(X,0,R))} \| -1| +$$
$$\sup_{|x|>R} |F_n(\gamma _{\xi ,n})(x)| \| \mu _{L(n)}|_{
(L(n)\setminus B(X,0,R))} \| .$$
Therefore, from
$\lim_{R \to \infty } \sup _{|x| >R} |F_n(\gamma _{\xi ,n})(x)|=0$
and from Lemma 2.3 the statement of Lemma 3.6 follows.
\par {\bf 3.7. Notes and definitions.}
Suppose $X$ is a locally convex space over a locally compact field
$\bf K$ with non-trivial non-Archimedean valuation
and $X^*$ is a topologically adjoint space.
The minimum $\sigma $-algebra with respect to which
the following family $ \{ v^* : v^*\in X^* \} $
is measurable is called a $\sigma $-algebra of cylinder sets.
For a $\bf K_s$-valued measure $\mu $ on $X$ a completion of
a linear space of characteristic functions $ \{ ch_U:
U\in Bco(X) \} $ in $L(X,\mu ,{\bf K_s})$ is denoted by 
$B_{\mu }(X)$.
Then $X$ is called a $RS$-space (or $KS$-space)
if on $X^*$ there exists a topology $\tau $ such that the 
continuity of each positive definite function $f: X^*\to \bf C$
(or $f: X^*\to \bf C_s$ with $\| f \|_{C^0} <\infty $
respectively) is necessary and sufficient for $f$ to be a characteristic 
functional of a non-negative measure
(a tight measure of finite norm correspondingly). Such topology
is called the $R$-Sazonov (or $K$-Sazonov) type topology.
The class of $RS$-spaces (and $KS$-spaces) contains
all separable locally convex spaces over $\bf K$.
For example, $l^{\infty }(\alpha ,{\bf K})=c_0(\alpha ,{\bf K})^*$,
where $\alpha $ is an ordinal \cite{roo}.
In particular we also write $c_0({\bf K}):=c_0(\omega _0,{\bf K})$ and
$l^{\infty }({\bf K}):=l^{\infty }(\omega _0,{\bf K})$,
where $\omega _0$ is the first countable ordinal.
\par Let $n_{\bf K}(l^{\infty },c_0)$ denotes the weakest topology on 
$l^{\infty }$ for which all functionals $p_x(y):=\sup_n|x_ny_n|$ are 
continuous, where $x=\sum_nx_ne_n\in c_0$ and $y=\sum_ny_ne^*_n
\in l^{\infty }$, $e_n$ is the standard base in $c_0$.
Such topology $n_{\bf K}(l^{\infty },c_0)$ is called the normal topology.
The induced topology on $c_0$ is denoted by $n_{\bf K}(c_0,c_0)$.
\par {\bf 3.8. Theorem.} {\it Let $f: l^{\infty }({\bf K})
\to \bf C$ (or $f: l^{\infty }({\bf K})\to \bf C_s$)
be a functional such that
\par $(i)$ $f$ is positive definite (or $f(0)=1$ and $ \| f \| _{C^0}\le 1$),
\par $(ii)$ $f$ is continuous in the normal topology
$n_{\bf K}(l^{\infty },c_0)$, then $f$ is the characteristic functional
of a probability measure on $c_0({\bf K})$.}
\par {\bf Proof.} If $\nu $ is the Haar measure on $\bf K^n$, then
on $Bco({\bf K^n})$ it takes values in $\bf Q$. Therefore, Lemma 4.1
\cite{mamaz} is transferable onto the case of $\bf K_s$-valued measures,
since ${\bf Q}\subset \bf K_s$. Therefore, analogously to Equation
$(4.1)$ of Lemma 4.2 \cite{mamaz} we have
$$(i)\mbox{ }P \{ |V_1|_{\bf K}<\epsilon ,...,
|V_n|_{\bf K}<\epsilon \} = \nu ^{-1}(B({\bf K^n},0,p^{-m}))
\int_{\bf K^n} f_V(y)ch_{B({\bf K^n},0,p^{-m})}(y)\nu (dy)$$
for measurable maps $V_j: (\Omega ,{\sf B},P)\to ({\bf K},Bco({\bf K}))$,
where $(\Omega ,{\sf B},P)$ is a probability space for a probability measure
$P$ with values in $\bf K_s$ on an algebra ${\sf B}$ of subsets of a set 
$\Omega $, $f_W$ is a characteristic function of $W=(V_1,...,V_n)$. 
To continue the proof we need the following statements.
\par {\bf 3.9. Lemma} {\it Let $f: c_0({\bf K})\to \bf C_s$ be a function
satisfying the following two conditions:
\par $(i)$ $|f(x)|\le 1$ for each $x\in c_0({\bf K})$,
\par $(ii)$ $f$ is continuous at zero in the topology 
$n_{\bf K}(c_0,c_0)$,  \\
then for each $\epsilon >0$ there exists $\lambda (\epsilon )\in c_0
({\bf K})$ such that $|1-f(x)|<p_{\lambda (\epsilon )}(x)+\epsilon $
for each $x\in c_0({\bf K}).$}
\par {\bf Proof.} In view of continuity for each $\epsilon >0$
there exists $y(\epsilon )\in c_0$ such that $|1-f(x)|<\epsilon $ 
if $p_{y(\epsilon )}<1$. Put $\lambda (\epsilon )=\pi _{\bf K}^{-1}
y(\epsilon )$, where $\pi _{\bf K}\in \bf K$ is such that 
$|\pi _{\bf K}|=p^{-1}$. If $x\in c_0$ is such that
$p_{\lambda (\epsilon )}(x)<p^{-1}$, then
$|1-f(x)|<\epsilon \le \epsilon + p_{\lambda (\epsilon )}(x).$
If $p_{\lambda (\epsilon )}(x)\ge p$, then $|1-f(x)|\le 2 \le p
<p_{\lambda (\epsilon )}(x)+\epsilon $.
\par {\bf 3.10. Lemma.} {\it Let $\{ V_n : n\in {\bf N} \} $ 
be a sequence of $\bf K$-valued random variables
for $P$ with values in $\bf K_s$.
If for each $\beta >0$ and $\epsilon >0$ there exists 
$N_{\epsilon }\in \bf N$ such that
$$(i)\mbox{ } \| P |_{
\{ \sup_{n\ge N_{\epsilon }} |V_n|_{\bf K}\le \beta  \} } \|
\ge 1-\epsilon (1+\beta ^{-1}),$$ then
$\lim_nV_n=0$ $P$-a.e. on $\Omega $.}
\par {\bf Proof} is quite analogous to that of Lemma 4.4 \cite{mamaz}
with substitution of $P$ on $\| P \| $.
\par {\bf 3.11. Proposition.} {\it Let $f: c_0({\bf K})\to \bf C_s$ 
be a function such that
\par $(i)$ $f(0)=1$ and $|f(x)|\le 1$ for each $x\in c_0$,
\par $(ii)$ $f(x)$ is continuous in the normal topology $n_{\bf K}(c_0,c_0)$.
Then there exists a probability measure $\mu $ on $c_0({\bf K})$
such that $f(x)={\hat \mu }(x)$ for each $x\in c_0$.}
\par {\bf Proof.} Consider functions $f_n(x_1,...,x_n):=f(
x_1e_1+...+x_ne_n),$ where $x=\sum_jx_je_j\in c_0$.
From Condition $(ii)$ and Proposition 3.1(2) \cite{mamaz}
it follows, that $f(x)$ is continuous in the norm topology.
From Chapters 7,9 \cite{roo} it follows, that there
exists a consistent family of tight measures $\mu _n$ on $\bf K^n$ such that
${\hat \mu }_n(x)=f_n(x)$ for each $x\in \bf K^n$.
In view of Theorem 3.1 there exists a probability space
$(\Omega ,{\sf B},P)$ with a $\bf K_s$-valued measure $P$ and a sequence 
of random variables $\{ V_n \} $ such that
$\mu _n(A)=P\{ \omega \in \Omega :$ 
$(V_1(\omega ),...,V_n(\omega ))\in A \} $ for each clopen subset 
$A$ in $\bf K^n$, consequently, $\lim_nV_n=0$ $P$-a.e. in $\Omega $.
In view of the preceding lemmas we have the following inequality:
$$|1- \| P|_{(|V_n|<\beta ,...,|V_{n+m}|<\beta )} \|
\le \| p_{\lambda (\epsilon )} (y_1e_n+...+y_me_{n+m}\| _{
L(B({\bf K^n},0,\beta ^{-1}),\nu ,{\bf K_s})}.$$
Since $\lim_kp_{\lambda (\epsilon )}(e_k)=0$, then there exists
$N\in \bf N$ such that $\sup_{k\ge N}p_{\lambda (\epsilon )}
(e_k)\le \epsilon $,
consequently, $ \| P |_{ \{ |V_N|<\beta ,...,|V_{N+m}||<\beta  \} } \|
\ge 1-\epsilon (1+\beta ^{-1})$. Due to Lemma 4.10
$ \| P |_{ \{ \lim_nV_n=0 \} } \| =1$.
Define a measurable mapping
$W$ from $\Omega $ into $c_0$ by the following formula:
$W(\omega ):=\sum_nV_n(\omega )e_n$ for each $\omega \in \Omega $,
then we also define a measure $\mu (A):=P \{ W^{-1}(B) \} $
for each $A\in Bco(X)$, hence $\mu $ is a probability measure on $c_0$.
In view of the Lebesgue convergence theorem (see Chapter 7 \cite{roo})
there exists ${\hat \mu }(x)=\lim_n{\hat \mu }_n(x_1e_1+...+x_ne_n)=f(x)$
for each $x\in c_0$.
\par {\bf Continuation of the proof of Theorem 3.8.}
Let $f: l^{\infty }({\bf K})\to \bf C_s$ satisfies
assumption of Theorem 3.8,  then by Proposition 3.11
there exists a probability measure $\mu $ on $c_0({\bf K})$
such that $f(x)={\hat \mu }(x)$ for each $x\in c_0({\bf K}).$
\par The case of the topological vector space $X$ over
$\bf K$ with $char({\bf K})>0$ and a real-valued measure
$\mu $ can be proved analogously to \cite{mamaz} 
due to \S 2.6 and \S \S 3.1-3.11.
\par {\bf 3.12. Theorem.} {\it Let $\mu $ be a probability measure
on $c_0({\bf K})$, then $\hat \mu $ is continuous in 
the normal topology $n_{\bf K}(l^{\infty },c_0)$ on $l^{\infty }$.}
\par {\bf Proof.} In view of Lemma 2.3 and Theorem 3.1
for each $\epsilon >0$ there exists $S(\epsilon )\in c_0$ such that
$\| \mu |_{L(0,S(\epsilon ))} \| \ge 1 -\epsilon $,
where $L(y,z):= \{ x\in c_0:$ $|x_n-y_n|\le |z_n|,$ for each 
$n\in {\bf N} \} $. Therefore, 
$$|1-{\hat \mu }(x)|
\le \epsilon +\| 2\pi \eta (\xi x) \|_{C^0(L(0,S(\epsilon )))} 
\| \mu |_{L(0,S(\epsilon ))} \| ,$$ 
hence there exists a constant
$C>0$ such that $|1-{\hat \mu }|\le \epsilon +Cp_{S(\epsilon )}(x)$.
\par {\bf 3.13. Corollary.} {\it The normal topology $n_{\bf K}
(l^{\infty },c_0)$ is the $R$-Sazonov (and $K$-Sazonov) type topology
on $l^{\infty }({\bf K})$.}
\par  {\bf 3.14. Theorem. Non-Archimedean analog of the
Minlos-Sazonov theorem.}
{\it For a separable Banach space 
$X$ over $\bf K$ the following two conditions are equivalent:
$$(I)\mbox{ }\theta : X \to {\bf T_s}
\mbox{ satisfies either conditions }2.6(3,4,5) \mbox{ or }2.6(3,6)
\mbox{ and }$$ for each $c>0$ there exists a compact operator
$S_c:X\to X$
such that either $|Re(\theta (y)-\theta (x))| <c$ 
or $|\theta (y)-\theta (x)|<c$ respectively for $|\tilde z (S_cz)|<1$;
$$(II)\mbox{ }\theta \mbox{ is a characteristic functional of a 
probability Radon measure } \mu $$ on $E$, where $\tilde z$ is an element
$z \in X
\hookrightarrow X^*$ considered as an element of $X^*$ under the
natural embedding associated with the standard base of 
$c_0(\omega _0,{\bf K})$, $z=x-y$, $x$
and $y$ are arbitrary elements of $X$.}
\par {\bf Proof.} $(II\to I)$. For a positive definite function $\theta $
generated by a probability measure $\mu $ in view of the inequality
$|\theta (y)-\theta (x)|^2 \le 2 \theta (0)(\theta (0)-Re(\theta
(y-x))$ (see Propositions IV.1.1(c)\cite{vah}) and using the
normalization of a measure $\mu $ by $1$
we consider the case $y=0$. For each $r>0$ we have: $|Re(\theta (0)-
\theta (x))|=\int_X (1-cos(2 \pi \eta (x(u)))) \mu (du)$ $\le
\int_{B(X,0,r)} 2 sin^2(\pi \eta (x(u))) \mu (du)$ $+2\int_{
X\setminus B(X,0,r)} $ $\mu (du)$ $\le 2 \pi ^2 \int_{B(X,0,r)}
\eta (x(u))^2 \mu (dx)$ $+2 \mu([x:$ $\| x\| >r])$. 
For $\theta $ generated by a $\bf K_s$-valued measure for each $r>0$
we have $| \theta (0)-\theta (x)|=|\int_X(1-exp[2\pi i \eta (x(u))])
\mu (du)|\le \| (1-exp[2 \pi
i \eta (x(u))])|_{B(X,0,r)} \| _{\mu }+2 \| \mu |_{(X
\setminus B(X,0,r))} \| $. In view of the
Radon property of the space $X$ and Lemma 2.5 for each $b>0$ and
$\delta >0$ there are a finite-dimensional over $\bf K$ subspace
$L$ in $X$ and a compact subset $W \subset X$ such that
$W \subset L^{\delta }$, $\| \mu |_{(X\setminus W)} \| <b$, 
hence $\| \mu |_{(X\setminus L^{\delta })} \| <b.$
\par We consider the following expression:
$$J(j,l):=2 \pi ^2 \int_{B(X,0,r)}
\eta (e_j(u)) \eta (e_l(u)) \mu (du),$$ 
where $(e_j)$ is the orthonormal basis
in $X$ which contains the orthonormal basis 
of $L=\bf K^n$, $n=dim_{\bf K}L$. Then we choose
sequences $b_j=p^{-j}$ and
$0< \delta _j <b_j$, subspaces $L_j$ and
$r=r_j$ such that $b_jr_j<1,$ $W_j \subset B(X,0,r_j)$, $0<r_j<r_{j+1}<
\infty $ for each $j \in \bf N$ and the orthonormal basis 
$(e_j)$ corresponding to the sequence
$L_j \subset L_{j+1} \subset ... \subset X$.
We get, due to finiteness of $n_j:=dim_{\bf K}L_j$, that
$\lim_{j+l \to \infty } J(j,l)=0$,
since $\| \mu |_{ \{ x: \| x\| >r_j \} } \| <b_j$,  
$\eta (x(u))=0$ for $x \in X \ominus
L_j$ with $\| x\| <b_j$, $u \in B(X,0,r_j)$. Then we define $g_{j,l}
:=\min \{ d:$ $d \in \Gamma _{\bf K}$  and  $d \ge |J(j,l)| \} $, evidently,
$g_{j,l} \le p|J(j,l)|$ and there are $\xi _{j,l} \in \bf K$
with $|\xi _{j,l}|_K=g_{j,l}$. Consequently, the family $(\xi _{j,l})$
determines a compact operator $S: X \to X$ with $\tilde e_j(Se_l)=
\xi _{j,l}t$
due to Theorem 1.2\cite{sch2}, where $t=const \in \bf K$, $ t\ne 0$. Therefore,
$|Re(\theta (0) -\theta (z))|<c/2 +|\tilde z(Sz)|<c$ 
for the real-valued measure $\mu $ and $|\theta (0)-\theta (z)|<c/2 + 
|{\tilde z}(Sz)|<c$ for the $\bf K_s$-valued measure,
if $|\tilde z(Sz)| <|t|c/2$. We choose $r$ such that
$\| \mu |_{(X\setminus B(X,0,r))} \| <c/2$ 
with $S$ corresponding to $(r_j:$ $j)$,
where $r_1=r$,
$L_1=L$, then we take $t \in \bf K$ with $|t|c=2$.
\par  $(I \to II)$. Without restriction of generality we may take
$\theta (0)=1$
after renormalization of non-trivial $\theta $. 
In view of Theorem 3.8 as in \S 2.6 
we construct using $\theta (z)$ a consistent family
of finite-dimensional distributions $\{ \mu _{L(n)} \} $
all with values either in $\bf R$ or $\bf K_s$ respectively. Let $m_{L(n)}$
be a real Haar measure 
on $L(n)$ which is considered as $\bf Q_p^a$ with
$a=dim_{\bf K}L(n)dim_{\bf Q_p}{\bf K}$, $m(B(L(n),0,1))=1$. In view of
Proposition 2.7 and Lemmas 2.8, 3.6: $\int_{L(n)} G_{\xi }(x) \mu
_{L(n)} (dx)$ $=\int _{L(n)} \gamma _{\xi ,n}(z) \theta (z)m_{L(n)}(dz)$,
consequently, 
$$1-\int_{L(n)} F_n(\gamma _{\xi ,n})(x) \mu _{L(n)} (dx)=
\int \gamma _{\xi ,n}(z) (1-\theta (z)) m_{L(n)} (dz)=:I_n(\xi ).$$
There exists an orthonormal basis 
in $X$ in which $S_c$ can be reduced to the following
form $S_c=SC\hat S_cE$ (see appendix), where $\hat S_c=diag(s_j:$
$j \in {\bf N})$ in the orthonormal basis 
$(f_j :j)$ in $X$ and $S$ transposes a finite
number of vectors in the orthonormal basis.
That is, $|\tilde z(\hat S_cz)|=\max_j |s_j|\times |z_j|^2$.
In the orthonormal basis $(e_j: j)$ adopted to
$(L(n):$ $n)$ we have $|\tilde z(S_cz)|=\max_{j,l
\in \bf N}(|s_{j,l}|\times |z_j| \times |z_l|)$, $\| S_c\| =\max_{j,l}
|s_{j,l}|$, where $S_c=(s_{j,l}:$ $j,l\in {\bf  N})$ in the 
orthonormal basis $(e_j)$, $r=const>0$. 
In addition, $p^{-1}|x|_K \le |x|_p \le p|x|_{\bf K}$ for each
$x \in \bf K$.
If $S_c$ is a compact operator such that  $|Re(\theta (y)-
\theta (x))|<c$ (or $|\theta (y)-\theta (x)|<c$)
for $|\tilde z(S_cz)|<1$, $z=x-y$, then either
$|Re(1-\theta (x))|< c+2|\tilde x(S_cx)|$ and 
$$I_n(\xi ) \le \int \gamma _{\xi ,n}
(z) [c+2|\tilde z(S_cz)|_K] m_{L(n)} (dz)\le c+b\| S_c\| /|\xi |^2,$$
$b=const $ is independent on $n$, $\xi $ and $S_c$, 
$$b:=p\times 
\sup_{|\xi |>r}
|\xi |^2 \int_{L(n)} \gamma _{\xi ,n}(z) |z|_p^2 m_{L(n)} (dz)< \infty $$
for the real-valued measures, or \\
$|1-\theta (x)||< \max (C, 2| {\tilde x}(S_cx)|)$ and
$ \| \gamma _{\xi ,n}(z)(1-\theta (z)) \| _{m_{L(n)}} \le $ \\ 
$\max ( \| \gamma _{\xi ,n}(z) \| _{m_{L(n)}}C,
2|(\gamma _{\xi ,n}(z)){\tilde z}(S_cz)|_{m_{L(n)}})\le $
$\max (C,b \| S_c \| /|\xi |^2)$, \\
where $b:=p\times \sup _{|\xi |>r}
(|\xi |^2 \| \gamma _{\xi ,n}(z) z^2 \| _{m_{L(n)}})<\infty $ for 
the $\bf K_s$-valued measures.
Due to the formula of changing variables in integrals (A.7\cite{sch1})
the following equality is valid: $J_n(\xi )=I_n(\xi )J_n(1)/
[I_n(1)|\xi |^2]$ for $|\xi |
\ne 0$, where $$J_n(\xi )=\int_{L(n)} \gamma _{\xi ,n}(z)|z|_p^2 m_{L(n)}
(dz).$$ Therefore, 
$$1-\int_X G_{\xi }(x) \mu _*(dx) \le c+b\| S_c
\| /|\xi |^2$$ 
for the real-valued measures
and $$|1- \| G_{\xi }(x) \| _{\mu _*}|\le \max (C,b \| S_c \|/|\xi |^2)$$
for the $\bf K_s$-valued measures.
Then taking the limit with $|\xi | \to \infty $ and then
with $c \to +0$ with the help of Lemma 3.6 we get the statement $(I\to II)$.
\section{The non-Archimedean analog of the Kakutani theorem.}
\par  {\bf 4.1. Definition.} Let on a completely regular space
$X$ with $ind (X)=0$ two non-zero real-valued (or $\bf K_s$-valued)
measures $\mu $ and $\nu $ are given.
Then $\nu $ is called absolutely continuous relative to $\mu $ if $\nu (A)=0$
for each $A \in Bf(X)$ with $\mu (A)=0$ (or there exists $f$ such that
$\nu (A)=\int _A
f(x)\mu (dx)$ for each $A\in Bco(X)$, where $f \in L(X,\mu ,{\bf K_s})$ 
respectively)
and it is denoted $\nu \ll \mu $. Measures $\nu $ and $\mu $ are
singular to each other if there is $F\in Bf(X)$ (or $F \in E$) with
$|\mu |(X\setminus F)=0$ and $|\nu |(F)=0$ (or
$ \| X\setminus F\| _{\mu }=0$ and $\| F\| _{\nu }=0$
respectively) and it is denoted
$\nu \perp \mu $. If $\nu \ll \mu $ and $\mu \ll \nu $ then they are called
equivalent, $\nu \sim \mu $.
\par  {\bf 4.2. Theorem.} {\it (A). Measures $\mu ^j: Bf(X)\to
\bf R$ (or $\mu ^j: E\to \bf K_s$), 
$j=1,2$, for a Banach space $X$ over $\bf K$ are orthogonal
$\mu ^1 \perp \mu ^2$
if and only if $\rho (x)=0 $ $(mod $ $\mu ^1)$ (or
$N_{{\mu ^1}}(x)N_{{\mu ^2}}(x)=0 $ for each $x \in X$ respectively).
\par (B). If for measures $\mu ^j: E \to \bf K_s$ on a
Banach space $X$ over $\bf K$
is satisfied $\rho (x)=0$ for each $x$ with $N_{{\mu ^1}}(x)>0$,
then $\mu ^1\perp \mu ^2$; the same is true for a completely regular
space $X$ with $ind(X)=0$
and $\rho (x)=\mu ^2(dx)/ \mu ^1(dx)=0$ for 
each $x$ with $N_{{\mu ^1}}(x)>0$.}
\par  {\bf Proof.} (A). In the case of real-valued $\mu ^j$
the proof differs only slightly from the proof of Theorem 2 \S 15\cite{sko}.
For $\mu ^j$ with values in $\bf K_s$ from Definition 4.1 it follows
that there exists $F \in E$ with $\| X\setminus F\|_{{\mu ^1}}=0$ and
$\|F\|_{{\mu ^2}}=0$. In view of Theorems 7.6 and 7.20\cite{sko}
the characteristic function $ch_F$ of the set $F$ belongs to
$L(\mu ^1)\cap L(\mu ^2)$
such that $N_{{\mu ^j}}(x)$ are semi-continuous from above,
$\| ch_F \|_{N_{\mu ^2}}
=0$, $\| ch_{X\setminus F} \|_{N_{{\mu ^1}}}=0$, consequently,
$N_{{\mu ^1}}(x)N_{{\mu ^2(x)}}=0$ for each $x \in X$.
\par  On the other hand, if $N_{{\mu ^1}}(x)N_{{\mu ^2}}(x)=0$ for each
$x$, then for $F:=[x \in X:$ $N_{\mu ^2}(x)=0]$ due to Theorem 7.2 \cite{roo}
$\|F\|_{{\mu ^2}} =\|ch_F\|_{N_{{\mu ^2}}}=0$. 
Moreover, in view of Theorem 7.6\cite{roo} $F=
\bigcap_{n=1}^{\infty } U_{s^{-n}}$, where $U_c:=[x \in X:$ $N_{\mu ^2}
(x)<c]$ are open in $X$, hence $ch_F \in L(\mu ^1)\cap L(\mu ^2)$
and $N_{\mu ^1}|_{(X\setminus F)}=0$, consequently, $\|X\setminus
F\|_{\mu ^1}=0$.
\par  (B). In view of Theorem 3.1(B) for each $A \in P_{L(n)}^{-1}
[E(L(n))]$ and $m>n$: $\int_A \rho _m(x)\mu ^1(dx)=\mu ^2(A)$, then from
$\lim_{n \to \infty }\| \rho (x)- \rho _n(P_{L(n)}x)\|
_{N_{\mu ^1}}=0$ and Conditions 2.1.(i-iii) on $\mu ^2$
Statement (B) follows.
\par  {\bf 4.3. Note.} For real-valued measures $\mu ^j$ on $Bf(X)$
for a Banach space $X$ over $\bf K$ (instead of a Hilbert space) 
using the above
given statements Theorems 3-6 and corollary in \S 15\cite{sko} may be
reformulated and proved.
The Radon-Nikodym theorem is not valid for  $\mu ^j$ with values in
$\bf K_s$, so not all theorems for real-valued measures may be
transferred onto this case. Therefore, the definition
of absolute continuity of measures was changed (see \S 4.1).
\par  {\bf 4.4. Theorem.} {\it Let measures $\mu ^j$ and $\nu ^j$
be with values in $\bf K_s$ on $Bf(X_j)$ for a 
Banach space $X_j$ over $\bf K$
and $\mu =
\mu ^1 \otimes \mu ^2$, $\nu =\nu ^1 \otimes \nu ^2$ on $X=X_1
\otimes X_2$, therefore, the statement $\nu \ll \mu $ is equivalent to
$\nu ^1 \ll \mu ^1$ and $\nu ^2 \ll \mu ^2$, moreover,
$\nu (dx)/ \mu (dx)=
(\nu ^1(P_1dx)/ \mu ^1 (P_1dx))(\nu ^2(P_2dx)/ \mu ^2(P_2dx))$, where
$P_j:X \to X_j$ are projectors.}
\par  {\bf Proof} follows from Theorem 7.15\cite{roo}
and modification of the proof of Theorem 5 \S 15\cite{sko}.
\par   {\bf 4.5. Theorem. The non-Archimedean analog of 
the Kakutani theorem.} {\it Let $X=\prod _{j=1}^{\infty } X_j$
be a product of completely regular spaces with
$ind(X_j)=0$ and probability measures $\mu ^j,$ $\nu ^j: E(X_j)\to \bf K_s$,
also let $\mu _j \ll \nu _j$ for each
$j$, $\nu =\bigotimes_{j=1}^{\infty }\nu _j$, $\mu = \bigotimes _{j=1}
^{\infty }\mu _j$ are measures on $E(X)$,
$\rho _j(x)=\mu _j(dx)/\nu _j(dx)$
are continuous by $x \in X_j$, $\prod_{j=1}^n \rho _j(x_j)=:t_n(x)$
converges uniformly on
$Af(X, \mu )$-compact subsets in $X$, $\beta _j
:=\| \rho _j(x)\|_{\phi _j}$, $\phi _j(x):=N_{\nu ^j}(x)$ on $X_j$.
If $\prod _{j=1}^{\infty } \beta _j$ converges in $(0, \infty )$
(or diverges to 0), then $\mu \ll \nu $ and $q_n(x)=\prod _{j=1}^n\rho _j
(x_j)$ converges in $L(X,\nu ,{\bf K_s})$ 
to $q(x)=\prod _{j=1}^{\infty }\rho _j(x_j)$
$=\mu (dx)/ \nu (dx)$ (or $\mu \perp \nu$ respectively), where $x_j
\in X_j$, $x \in X$.}
\par  {\bf Proof.} The countable additivity of $\nu $ and $\mu $
follows from Theorem 3.1. Then $\beta _j=\| \rho _j\| _{\phi _j}
\le \| \rho _j\| _{N_{\nu _j}}=\| X \| _{\mu _j}=1$, since $N_{\nu _j}
\le 1$ for each $x \in X_j$, hence $\prod _{j=1}^{\infty } \beta _j$
can not be divergent to $\infty $. If this product diverges to $0$
then there exists a sequence $\epsilon _b:=\prod_{j=n(b)}^{m(b)}
\beta _j$ for which the series converges $\sum _{b=1}^{\infty } \epsilon _b
<\infty $,
where $n(b) \le m(b)$. For $A_b:=[x:$ $(\prod_{j=n(b)}^{m(b)} \rho _j(x_j)
) \ge 1]$ there are estimates $\| A_b \| _{\nu } \le $ $\sup_{x \in A_b}
[\prod_{j=n(b)}^{m(b)}| \rho _j(x_j)| \phi _j(x_j)] \le \epsilon _b$,
consequently, $\| A \| _{\nu }=0$ for $A=\lim \sup (A_b:$ $b \to \infty )$,
since $0< \sum_{b=1}^{\infty } \epsilon _b < \infty $.
\par  For $B_b:=X\setminus A_b$ we have: $\| B_b \| _{\mu } \le $ $[\sup_{x
\in B_b} \prod _{j=n(b)}^{m(b)} |1/ \rho _j(x_j)| \psi (x_j)]$
$=[\prod_{j=n(b)}^{m(b)} \| \rho _j(x_j) \| _{\phi _j}]=
\epsilon _b$, where $\psi _j(x)=N_{\mu _j}(x)$, since $\mu _j(dx_j)
= \rho _j(x_j) \nu _j(dx_j)$ and $N_{\mu _j}(x)=|\rho _j(x_j)|
N_{\nu _j}(x)$ due to continuity of $\rho _j (x_j)$ (for $\rho _j
(x_j)=0$ we set $|1/ \rho _j(x_j)| \psi _j(x_j)=0$, because
$\psi _j(x_j)=0$ for such $x_j$), consequently, $\| \lim \sup (B_b:$
$ b \to \infty )\| _{\mu }=0$ and $\| A \| _{\mu } \ge \| \lim \inf
(A_b:$ $b \to \infty) \| _{\mu }=1$. This means that $\mu \perp \nu $.
\par  Suppose that $\prod_{j=1}^{\infty } \beta _j$ converges to
$0<\beta < \infty $, then $\beta \le 1$ (see above). 
Therefore from the Lebesgue
Theorem 7.F\cite{roo} it follows that $t_n(x)$ converges in 
$L(X,\mu ,{\bf K_s})$, since
$|t_n(x)| \le 1$ for each $x$ and $n$, at the same time
each $t_n(x)$ converges uniformly on compact subsets in the topology
generated by $Af(X, \mu )$.
Then for each bounded continuous cylinder function
$f: X\to \bf K_s$ we have 
$$\int_X f(x) \mu (dx)=\int_X f(x_1,...,x_n)
t_n(x) \otimes _{j=1}^n \nu _j(dx_j)=\lim _{j \to \infty }
\int_X f(x)t_n(x)\nu (dx)=\int_X \rho (x) \nu (dx).$$ Approximating
arbitrary $h \in L(X,\mu ,{\bf K_s})$ by such $f$ we get the equality
$$\int_X h(x) \mu (dx)=\int_X h(x) \rho (x) \nu (dx),$$ consequently,
$\rho (x)=\mu (dx)/\nu (dx)$.
\section{Appendix.}
Suppose $X=c_0(\omega _0,{\bf K})$ is a Banach space over $\bf K$
and  $I$ is a unit operator on $X$.
If $A$ is an operator on $X$, then in some basis of $X$ we have
an infinite matrix $(A_{i,j})_{i,j \in \bf N}$, so we can consider
its transposed matrix $A^t$. If in some basis the following equality 
is satisfied $A^t=A$, then $A$ is called symmetric.
\par {\bf A.1. Lemma.} {\it Let $A: X\to X$
be a linear invertible operator with a compact operator $(A-I)$.
Then there exist an orthonormal basis
$(e_j:$ $j\in {\bf N})$ in $X$, invertible linear
operators $C, E, D: X\to X$ with compact opeartors $(C-I)$, $(E-I)$, $(D-I)$
such that $A=SCDE$, $D$ is diagonal, $C$ is lower triangular and $E$ is
upper triangular, $S$ is an operator transposing a 
finite number of vectors from an orthonormal basis 
in $X$. Moreover, there exists $n \in \bf N$ and invertible
linear operators $A', A": X\to X$ with compact operators $(A'-I)$,
$(A"-I)$ and $(A'_{i,j}-\delta _{i,j}=0)$
for $i$ or $j>n$, $A"$ is an isometry and there exist their determinants
$det(A')det(A")=det(A)$, $|det(A")|_{\bf K}=1$, $det(D)=det(A)$.
If in addition $A$ is symmetric, then $C^t=E$ and $S=I$.}
\par  {\bf Proof.} In view of Lemma 2.2\cite{sch2} for each $c>0$
there exists the following decomposition $X=Y\oplus Z$ into
$\bf K$-linear spaces such that $\| (A-I)|_Z \| <c$, where
$dim_{\bf K}Y=m< \aleph _0$. In the orthonormal basis 
$(e_j:$
$j)$ for which $sp_{\bf K}(e_1,...,e_m)=Y$ for $c \le 1/p$ we get
$A=A'A"$ with $(A-I)|_Z=0$, $|A"_{i,j}-\delta _{i,j}| \le c$ for each
$i, j$ such that $(A'_{i,j}-\delta _{i,j})=0$ for $i$ or $j>n$, where
$n \ge m$
is chosen such that $|A_{i,j}-\delta _{i,j}|\le c^2$ for $i>n$ and
$j=1,...,m$,
$A_{i,j}:=e_i^*(Ae_j)$, $e_i^*$ are vectors $e_i$
considered as linear continuous functionals $e_i^*
\in X^*$. Indeed, $(A_{i,j}:$  $i \in {\bf N})=Ae_j\in X$
and $\lim_{i \to \infty }A_{i,j}=0$ for each $j$. From the form of
$A"$ it follows that $\| A"e_j-e_j\| \le 1/p$ for each $j$,
consequently, $\| A"x\|=\| x\|$ for each $x \in X$. Since
$A"=(A')^{-1}A$, $(A-I)$ and $(A'-I)$ being compact, hence $(A"-I)$
is compact together with $(A^{-1}-I)$, $((A')^{-1}-I)$ and
$((A")^{-1}-I)$. Moreover, there exists $\lim_{k \to \infty }
det(A)_k =det(A)$ $=\lim_k
det((A')_k(A")_k)$ $=\lim_kdet(A')_kdet(A")_k$ $=det(A')det(A")$,
where $(A)_k:=(A_{i,j}:$ $i,j \le k)$. This follows from the
decompositions
$X=Y_k\oplus Z_k$ for $c=c(k) \to 0$ whilst $k \to \infty $. This
means that for each
$c(k)=p^{-k}$ there exists $n(k)$ such that $|A_{i,j}-
\delta _{i,j} |<c(k)$, $|A'_{i,j}-\delta _{i,j}|<c(k)$ and $|A"_{i,j}-
\delta _{i,j}|<c(k)$ for each $i$ or $j >n(k)$, consequently,
$|A{{1...n(k)i_1...i_q}\choose {1...n(k)j_1...j_q}} -A{{1...n(k)}
\choose {1...n(k)}}\delta _{i_1,j_1}...\delta _{i_q,j_q}| <c(k)$,
where $A{{i_1...i_r}\choose {j_1...j_r}}$ is a minor corresponding to
rows $i_1,...,i_r$ and columns $j_1,...,j_r$ for $r, q \in \bf N$. From
the ultrametric inequality it follows that $|det(A")-1| \le 1/p$,
hence $|det(A")|_{\bf K}=1$, $det(A")_k \ne 0$ for each $k$, $det(A')_k
=det(A')_n$ for each $k \ge n$. Using the decomposition of
$det(A')_n$
by the last row (analogously by the column) we get
$A'_{n,j} \ne 0$ and a minor
$A'{{1...n-1}\choose {1...j-1,j+1,n}} \ne 0$. Permuting the columns
$j$ and $n$ (or rows) we get as a result a matrix $(\bar A')_n$ with $\bar
A'{{1...n-1}\choose {1...n-1}} \ne 0$. Therefore, by the
denumeration of the basic vectors we get $A'{{1...k}\choose {1...k}}
\ne 0$ for each $k=1,...,n$, since $|det(\bar A')_n|=|det(A')_n|$.
\par   Therefore, there exists the orthonormal basis 
$(e_j:$ $j)$ such that $A{{1...j}\choose
{1...j}} \ne 0$ for each $j$ and $\lim_j A{{1...j}\choose {1..j}}
=det(A)\ne 0$. Applying to $(A)_j$ the Gaussian decomposition and using
compactness of $A-I$ due to formula (44) in \S II.4\cite{gan}, which is
valid in the case of $\bf K$ also, we get $D=diag(D_j:$
$j\in {\bf N})$, $D_j=A{{1...j}\choose {1...j}}/A{{1...j-1}\choose
{1...j-1}}$; $C_{g,k}=A{{1,...,k-1,g}\choose {1,...,k-1,k}}/A{{1...k}
\choose {1...k}}$; $E_{k,g}=A{{1,...,k-1,k}\choose {1,...,k-1,g}}/
A{{1...k}\choose {1...k}}$ for $g=k+1, k+2,...$, $k \in \bf N$.
Therefore, $(C-I)$, $(D-I)$, $(E-I)$ are the compact operators,
$C_{i,j}$, $D_j$,
$E_{i,j} \in \bf K$ for each $i, j$. Particularly, for $A^t=A$ ($A^t$
denotes the transposed matrix for $A$) we get $E_{k,g}=C_{g,k}$.

\par Address: Theoretical Department,
\par Institute of General Physics,
\par Russian Academy of Sciences,
\par Str. Vavilov 38, Moscow, 117942, Russia
\end{document}